\documentclass[9pt,shortpaper,twoside,web]{ieeecolor}
\usepackage{generic}
\pdfoutput=1
\usepackage{cite}
\usepackage{amsmath,amssymb,amsfonts}
\usepackage{algorithmic}
\usepackage[ruled,linesnumbered]{algorithm2e}
\usepackage{graphicx}
\usepackage{textcomp}
\usepackage{soul}
\def\BibTeX{{\rm B\kern-.05em{\sc i\kern-.025em b}\kern-.08em
    T\kern-.1667em\lower.7ex\hbox{E}\kern-.125emX}}
\begin{document}
\title{Accelerated Distributed Aggregative Optimization}
\author{Jiaxu Liu,  Song Chen , Shengze Cai,  Chao Xu, \IEEEmembership{Member, IEEE}
\thanks{This work was supported by the  Science and Technology Innovation 2030 New Generation Artificial Intelligence Major Project  NO.2018AAA0100902, the National Key Research and Development Program of China under Grant 2019YFB1705800, the National Natural Science Foundation of China number 61973270 and the Zhejiang Provincial Natural Science Foundation of China number LY21F030003.}
\thanks{Jiaxu Liu is with the School of Mathematical Sciences, Zhejiang University, Hangzhou, Zhejiang 310027, China. (e-mail: jiaxuliu@zju.edu.cn). }
\thanks{Song Chen is with the School of Mathematical Sciences, Zhejiang University, Hangzhou, Zhejiang 310027, China. (e-mail: math\_cs@zju.edu.cn).}
\thanks{Shengze Cai is with the State Key Laboratory of Industrial Control Technology, Institute of Cyber-Systems and Control, Zhejiang University, Hangzhou, Zhejiang 310027, China. (e-mail: shengze\_cai@zju.edu.cn).}
\thanks{Chao Xu is with the State Key Laboratory of Industrial Control Technology, Institute of Cyber-Systems and Control, Zhejiang University, Hangzhou, Zhejiang 310027, China, and also with Huzhou Institute of Zhejiang University, Huzhou, Zhejiang 313000, China.  (e-mail: cxu@zju.edu.cn).}
}

\maketitle

\begin{abstract}
In this paper, we investigate   a distributed  aggregative optimization problem in a network, where each agent has its own local cost function which depends  not only on the local state variable but also on an aggregated function of state variables from all agents. To accelerate the optimization process, we combine heavy ball and Nesterov's accelerated methods with distributed aggregative gradient tracking, and propose two novel algorithms named DAGT-HB and DAGT-NES for solving the distributed  aggregative optimization problem. We analyse that the DAGT-HB and DAGT-NES algorithms can converge to an optimal solution at a global $\mathbf{R}-$linear convergence rate when the objective function is smooth and strongly convex, and when the parameters (e.g., step size and momentum coefficients) are selected within certain ranges. A numerical
experiment on the optimal placement problem is given to verify the effectiveness and superiority of our proposed algorithms. 
\end{abstract}


\section{Introduction}
\label{sec:introduction}
Distributed optimization plays a critical role in machine learning, especially in the scenarios where the dataset or model is too large to fit in a single machine. By distributing the dataset or model over parallel machines and applying distributed optimization, one can achieve more efficient computation and scalable machine learning \cite{dean2012large,Barbarossa2014heterogeneous,li2014communication,Predd2009ACT}. Moreover,  distributed optimization has received  immense attention in the field of control because of its wide applications including formation control\cite{WANG2014nonholonomic},  sensor
networks \cite{Zhu2013Sensor}, resource allocation\cite{Deng2018DistributedCA}  and so on. Generally speaking, in a distributed optimization problem, each agent can only have access to its own information and communicate data with its neighbors to minimize a global cost function cooperatively. These agents can work in parallel through communication and collaboration, making the optimization process faster and able to handle larger datasets.

In order to solve the distributed optimization problem, various remarkable algorithms have been developed. Based on a consensus scheme, some known distributed optimization algorithms  are usually based on gradient descent such as the distributed subgradient descent \cite{Nedic2009Subgradient}, the distributed
dual averaging gradient algorithm\cite{Duchi2010DualAF}, push-sum distributed algorithm\cite{Tsianos2012PushSumDD} and so on. However, the above algorithms all have the disadvantage of slow convergence due to the use of a gradual vanishing step size in the algorithm design. With a constant step size, the convergence rate of the distributed gradient-based algorithm can be improved but they can only converge to 
a small neighborhood of the optimal solution due
to the use of the local gradient in each agent\cite{yuan2016convergence}.  To deal with the issues of slow convergence rate and  non-optimal solution, the gradient tracking strategy is merged with
the distributed convex optimization algorithms where an
estimate of the global average gradient is used to replace the local gradient in each agent\cite{xu2015augmented,Nedi2016AchievingGC}. In \cite{shi2015extra}, an algorithm called EXTRA was proposed to achieve geometric convergence to the global optimal solution by introducing a cumulative correction term.

Furthermore, in order to speed up the convergence rate, some acceleration algorithms such as heavy ball \cite{polyak1964some} and  Nesterov's accelerated method\cite{nesterov1983method}  can be applied to distributed optimization. In \cite{jakovetic2014fast}, the distributed Nesterov gradient method (D-NG) was proposed, which can improve the convergence rate to $O(\frac{\log k}{k})$. In \cite{xin2019distributed}, two
distributed Nesterov gradient methods over strongly connected and directed networks were proposed by extending $AB$~\cite{Ran2018Geometric} with Nesterov’s momentum. Subsequently,  in \cite{Ran2020Generalization}, the authors combined $AB$ with heavy ball and proposed a distributed heavy-ball algorithm that is named as $ABm$. It is also proven that $ABm$ has a global $R-$linear rate when the step-size and momentum parameters are positive and sufficiently small. Qu and Li \cite{Qu2020Accelerated} further combined the gradient tracking with distributed Nesterov's gradient descent method, resulting in two accelerated distributed Nesterov's methods termed as Acc-DNGD-SC and Acc-DNGD-NSC. Furthermore, a distributed Nesterov-like gradient tracking algorithm, called D-DNGT, which includes the gradient
tracking into the distributed Nesterov method with momentum
terms and employs nonuniform step sizes, was introduced in~\cite{lu2020nesterov}. These algorithms achieve linear convergence rates for smooth and strongly convex objective functions, and they significantly improve the convergence speed compared to algorithms without introducing momentum terms.

In the distributed optimization problems mentioned above, the main setup that has emerged is named consensus optimization. However, in some other scenarios, such as the multi-agent formation control problem, multi-robot surveillance scenario\cite{Carnevale2022AggregativeFO} and aggregative game \cite{koshal2016distributed}, the objective function of each agent is not only dependent to its local state but also determined by other agents’ variables through an aggregative variable.
Such an optimization problem is called distributed aggregative optimization in \cite{Li2022Aggregative}, where a distributed aggregative gradient tracking (DAGT) method was also proposed to handle this problem. 
Following \cite{Li2022Aggregative}, Li et.al \cite{Li2022Online} considered the online convex optimization with an aggregative variable to solve time-varying cost functions. Subsequently, the online distributed aggregative optimization problem with constraints was investigated in \cite{Carnevale2022Coordination}. Chen and Liang \cite{Chen2022DistributedAO} considered finite bits communication and proposed a  novel distributed quantized algorithm  called D-QAGT. In \cite{Wang2022DistributedPA}, the distributed convex aggregative optimization was further combined with Frank-Wolfe algorithm to solve the aggregative optimization problem over
time-varying communication graphs.

To the best of our knowledge,  little work has been
done to accelerate the convergence rate of distributed aggregative optimization. Therefore, we propose two novel accelerated algorithms for distributed aggregative optimization problem in this paper. The main contributions of our work can be summarized as follows. 
\begin{enumerate}
    \item  We combine accelerated algorithms with Distributed Aggregative Gradient Tracking (DAGT)\cite{Li2022Aggregative}, resulting in two novel algorithms called DAGT-HB and DAGT-NES.
    \item  We theoretically show that the algorithms DAGT-HB  and DAGT-NES can converge to an optimal solution at a global $R-$linear convergence rate when the objective function is smooth and strongly convex. Moreover, the proper ranges for selecting parameters (e.g., step size and momentum term) are provided. 
    \item  The numerical simulation verifies the effectiveness and superiority of the proposed DAGT-HB and DAGT-NES, supporting our theoretical findings in this paper.
\end{enumerate}

The rest of this paper is organized as follows. In Section 2, the  basic notations, the basic definition of graph and the distributed aggregative optimization
problem are presented. The DAGT-HB and DAGT-NES are introduced and their convergence rates are also analysed in Section 3 and Section 4, respectively. A numerical example is given in Section 5 to validate the proposed algorithms. Eventually,  Section 6 concludes the paper.

\section{Preliminaries And Problem Setup}
\subsection{Basic Notations}
The set of real and positive real numbers are denoted by $\mathbf{R}$ and $\mathbf{R}^{\mathbf{+}}$, the set of $n$-dimensional column vectors is denoted by $\mathbf{R}^n$. Let $\mathbf{1}_n$ and $\mathbf{0}_n$ represent the column vectors of $n$ ones and zeros, respectively. $I_n$ denotes the $n$-dimensional identity matrix. Let $||\cdot||$ and $x^\top$ be the standard Euclidean norm (or induced matrix norm) and the transpose of $x\in \mathbf{R}^n $.
We use pointwise order $>$ for any vectors $a=(a_1, a_2, \ldots, a_n)^\top$ and $b=(b_1, b_2, \ldots, b_n)^\top \in \mathbf{R}^n$, i.e., $a> b \iff a_i>b_i$,  $i=1, 2, \ldots, n$. For vectors $x_1, \cdots, x_N \in \mathbf{R}^n$, we use the notation $x=\operatorname{col}\left(x_1, \cdots\right.$, $\left.x_N\right)=\left(x_1^{\top}, \cdots, x_N^{\top} \right)^{\top}$ to denote a new stacked vector. Also, we  use
blkdiag$(A_1, A_2, \ldots, A_N)$ to represent the block diagonal matrix
where the $i$-th diagonal block is given by the matrix $A_i \in \mathbf{R}^{m_i\times n_i}$, $i=1, 2, \ldots, n$. The Kronecker product of arbitrary matrices $A \in \mathbf{R}^{m \times n}$ and $B \in \mathbf{R}^{p \times q}$ is defined as $A \otimes B \in$ $\mathbf{R}^{m p \times n q}$. Let $\rho(A)$ be the spectral radius of a square matrix $A$. $A>0$ represents that $A$ is positive, that is, every entry of $A$ is greater than 0. Let $K_N=\frac{1}{N}1_N1^\top_N$. Besides, $\nabla f(\cdot)$ is the gradient of a differentiable function $f(\cdot)$.

\subsection{Graph Theory}
Here, we provide some basic definitions of graph theory. Let $\mathcal{G}=(\mathcal{V}, \mathcal{E}, A)$ denote a weighted undirected graph with a finite vertex set $\mathcal{V}=\{1, \ldots, N\}$, an edge set $\mathcal{E} \subseteq \mathcal{V} \times \mathcal{V}$, and a weighted adjacency matrix $A=\left[a_{i j}\right] \in \mathbb{R}^{N \times N}$ with $a_{i j}>0$ if $(j, i) \in \mathcal{E}$ and $a_{i j}=0$ otherwise. $\mathcal{N}_i=\{j:(j, i) \in \mathcal{E}\}$ denotes the set of neighbors of agent $i$ and $d_i=\sum_{j=1}^N a_{i j}$ denotes the weighted degree of vertex $i$. The
graph $\mathcal{G}$ is called connected if for any $i, j \in  \mathcal{V}$ there exists a  path from
$i$ to $j$ . The Laplacian matrix of graph $\mathcal{G}$ is $L=D-A$ with $D=\operatorname{diag}\left\{d_1, \ldots, d_N\right\}$. The real eigenvalues of $L$ are denoted by $\lambda_1,...,\lambda_N$ with $\lambda_i \leq \lambda_{i+1}$, $i=1, \ldots, N-1$. The undirected graph $\mathcal{G}$ is connected if and only if $\lambda_2>0$.

\subsection{Problem Formulations}
In this paper, we consider the distributed aggregative
optimization problem that can be written as:
\begin{align}
    \min _{x \in \mathbb{R}^n} f(x)  =\sum_{i=1}^N f_i\left(x_i, u(x)\right) , ~~~  u(x)  =\frac{\sum_{i=1}^N \phi_i\left(x_i\right)}{N},
\end{align}
where $x=\operatorname{col}\left(x_1, \ldots, x_N\right)$ is the global state variable with $x_i \in$ $\mathbb{R}^{n_i}, n:=\sum_{i=1}^N n_i$, and $f_i: \mathbb{R}^{n_i} \rightarrow \mathbb{R}$ is the local objective function. In problem (1), $u(x)$ is an aggregative variable that can have access to information of all agents and the function $\phi_i$: $\mathbb{R}^{n_i} \rightarrow \mathbb{R}^d$ is only accessible to agent $i$. Moreover, each agent $i$ only knows the information of state variable $x_i$ and can not obtain the information of other state variables. And each agent can only privately access the information on $f_i$. The purpose of this paper is to design a distributed optimization algorithm to obtain the optimal state variable for problem (1).

For simplicity in the following analysis, let $\nabla_1 f_i\left(x_i, u(x)\right)$ and $\nabla_2 f_i\left(x_i, u(x)\right)$ denote $\nabla_{x_i} f_i\left(x_i, u(x)\right)$ and $\nabla_u f_i\left(x_i, u(x)\right)$, respectively, for all $i\in \{1,2,\cdots,N\}$. For $x=$ col$\left(x_1, \ldots, x_N\right) \in \mathbb{R}^{n}$ and $y=$ col$\left(y_1, \ldots, y_N\right) \in \mathbb{R}^{N d}$, we define $f(x, y)=\sum_{i=1}^N f_i\left(x_i, y_i\right)$, $\nabla_1 f(x, y)=$ col$\left(\nabla_1 f_1\left(x_1, y_1\right), \ldots, \nabla_1 f_N\left(x_N, y_N\right)\right)$ and $\nabla_2 f(x, y)=$ col$\left(\nabla_2 f_1\left(x_1,y_1\right), \ldots, \nabla_2 f_N\left(x_N, y_N\right)\right)$. Also, we denote   $\nabla \phi(x)=$ blkdiag$(\nabla\phi_1\left(x_1\right), \ldots, \nabla\phi_N\left(x_N\right))\in \mathbf{R}^{n\times Nd}$. Next, for a differentiable function $h(x)=\operatorname{col}\left(h_1(x), \ldots, h_m(x)\right)$, with $h_i$ being a real-valued function, let us denote $\nabla h(x)=\left(\nabla h_1(x), \ldots, \nabla h_m(x)\right)$ with $\nabla h_i$ being the gradient (a column vector) of $h_i, i=1,2,\ldots,m$. 

To facilitate the subsequent analysis, it is necessary to make some definitions and general assumptions. Firstly, we give the definitions of $L$-smooth function, $m$-strongly convex function and $\mathbf{R}$-linear convergence respectively.

\textbf{Definition $1$ (smoothness)}: A differentiable function $f: \mathbb{R}^n \rightarrow \mathbb{R}$ is $L$-smooth if for all $x, y \in \mathbb{R}^n$
\begin{equation}
    \|\nabla f(x)-\nabla f(y)\| \leq L\|x-y\|.
\end{equation}
Following \cite{nesterov2003introductory}, it is equivalent to that for all $x, y \in \mathbb{R}^n$
\begin{equation}
    f(y) \leq f(x)+\nabla f(x)^{\top}(y-x)+\frac{L}{2}\|y-x\|^2.
\end{equation}

\textbf{Definition $2$ (strong convexity)}: A differentiable function $f: \mathbb{R}^n \rightarrow \mathbb{R}$ is $m$-strongly convex  if for all $x, y \in \mathbb{R}^n$
\begin{equation}
    \quad m\|x-y\|^2 \leq(x-y)^{\top}(\nabla f(x)-\nabla f(y)). 
\end{equation}
 It also implies that for all $x, y \in \mathbb{R}^n$
\begin{equation}
    f(x)+\nabla f(x)^{\top}(y-x)+\frac{m}{2}\|y-x\|^2 \leq f(y).  
\end{equation}

\textbf{Definition $3$  }\cite{ortega2000iterative}: A sequence $\{x_k\}$  is said to converge $\mathbf{R}$-linearly to $x_*$
 with rate $\rho \in (0,1)$ if there is a constant $c>0$ such that
\begin{equation}
    ||x_k-x_*||\leq c\rho^k\quad \forall {k \in N_+}.
\end{equation}

Next, we make some common assumptions about the graph $\mathcal{G}$  and objective function $f(x)$.

\textbf{Assumption $1$}: The graph $\mathcal{G}$ is connected and $A$ is doubly stochastic, that is, $\sum_{i=1}^{N}a_{ij}=1$ and  $\sum_{j=1}^{N}a_{ij}=1$ for all $i,j=1,\ldots,N$.

\textbf{Assumption $2$}: The global objective function $f$ is $L_1$-smooth and $m$-strongly convex. This means that $\nabla_1 f(x, y)+\nabla \phi(x) \mathbf{1}_N \otimes \frac{1}{N} \sum_{i=1}^N \nabla_2 f_i\left(x_i, y_i\right)$ is  $L_1$-Lipschitz continuous.

\textbf{Assumption $3$}: $\nabla_2 f(x, y)$ is $L_2$-Lipschitz continuous.

\textbf{Assumption $4$}: For all $i=1,2,\ldots,N$, the aggregation function $\phi_i$ is differentiable and $L_3$-Lipschitz continuous.

To perform  the following analysis, several key lemmas are listed below.

\textbf{Lemma $1$} \cite{Li2022Aggregative}: Let $F: \mathbb{R}^n \rightarrow \mathbb{R}$ be $\mu$-strongly convex and $L$-smooth. Then, $\|x-\alpha \nabla F(x)-(y-\alpha \nabla F(y))\| \leq(1-m \alpha)\|x-y\|$ for all $x, y \in \mathbb{R}^n$, where $\alpha \in(0,1 / L]$.

\textbf{Lemma $2$}  \cite{horn2012matrix}: Under Assumption 1,  for the adjacency
matrix $A$, the following properties hold:
\begin{enumerate}
    \item $\mathcal{A} \mathcal{K}=\mathcal{K} \mathcal{A}=\mathcal{K}$, where $\mathcal{A}=A \otimes I_d$, $\mathcal{K}=K_N \otimes I_d$.
    \item $\|\mathcal{A} x-\mathcal{K} x\| \leq \rho\|x-\mathcal{K} x\|$ for any $x \in$ $\mathbb{R}^{N d}$ and $\rho =\|A-K_N\|<1$.
    \item $||A-I_N|| \leq 2$.
\end{enumerate}

\textbf{Lemma $3$} \cite{horn2012matrix}: Let $W\in \mathbb{R}^{n\times n} $ be nonnegative and $x\in \mathbb{R}^{n}$ be positive. If $Wx<\lambda x$ with $\lambda>0$, then $\rho(W)<\lambda$.

\textbf{Lemma $4$}\cite{zheng2010extended}: Let $\mathbb{R}$ be the real number field, and $H\left(z\right)$ denote the degree $n$ real coefficient polynomial
$$
H\left(z\right)=a_n z^n+a_{n-1} z^{n-1}+\cdots+a_2 z^2+a_1 z+a_0,
$$
 where $a_0, a_1, \ldots, a_n \in \mathbb{R}$.
The Jury matrix of $H\left(z\right)$ can be written as

\begin{tabular}{c|lllllll} 
& $z^0$ & $z^1$ & $z^2$ & $\cdots$ & $z^{n-2}$ & $z^{n-1}$ & $z^n$ \\
\hline 1 & $a_0$ & $a_1$ & $a_2$ & $\cdots$  & $a_{n-2}$ & $a_{n-1}$ & $a_n$ \\
2 & $a_n$ & $a_{n-1}$ & $a_{n-2}$ & $\cdots$ & $a_2$ & $a_1$ & $a_0$ \\
3 & $b_0$ & $b_1$ & $b_2$ & $\cdots$  & $b_{n-2}$ & $b_{n-1}$ & \\
4 & $b_{n-1}$ & $b_{n-2}$ & $b_{n-3}$ & $\cdots$  & $b_1$ & $b_0$ & \\
5 & $c_0$ & $c_1$ & $c_2$ & $\cdots$  & $c_{n-2}$ & \\
6 & $c_{n-2}$ & $c_{n-3}$ & $c_{n-4}$ & $\cdots$  & $c_0$ & \\
$\vdots$ & $\vdots$ & $\vdots$ & $\vdots$ & & & &  \\
$2 n-5$ & $l_0$ & $l_1$ & $l_2$ & $l_3$ & & &  \\
$2 n-4$ & $l_3$ & $l_2$ & $l_1$ & $l_0$ & & &  \\
$2 n-3$ & $m_0$ & $m_1$ & $m_2$ & & & &  \\
\end{tabular}
where
$$
\begin{aligned}
& b_i=\left|\begin{array}{ll}
a_0 & a_{n-i} \\
a_n & a_i
\end{array}\right|, \quad i=0,1,2, \ldots, n-1,\\ & c_j=\left|\begin{array}{ll}
b_0 & b_{n-j-1} \\
b_{n-1} & b_j
\end{array}\right|, \quad j=0,1,2, \ldots, n-2, \\
 & \cdots \\&  m_0=\left|\begin{array}{ll}
l_0 & l_3 \\
l_3 & l_0
\end{array}\right|, \quad m_1=\left|\begin{array}{cc}
l_0 & l_2 \\
l_3 & l_1
\end{array}\right|, \quad m_2=\left|\begin{array}{cc}
l_0 & l_1 \\
l_3 & l_2
\end{array}\right| . 
&
\end{aligned}
$$
All the modulu of roots of a real coefficient polynomial $H\left(z \right)\left(n \geqslant 3, a_n>\right.$ 0) are less than 1 if and only if the following four conditions hold:
\begin{enumerate}
    \item $H\left(1 \right)>0$;
    \item $(-1)^n H\left(-1 \right)>0$;
    \item $\left|a_0\right|<a_n$;
    \item $\left|b_0\right|>\left|b_{n-1}\right|, \left|c_0\right|>\left|c_{n-2}\right|, \ldots, \left|l_0\right|>\left|l_3\right|, \left|m_0\right|>\left|m_2\right|$.
\end{enumerate}

\section{DAGT-HB}

In order to solve problem (1), we combine distributed aggregative gradient tracking (DAGT)\cite{Li2022Aggregative} with heavy ball method and propose the following DAGT-HB algorithm:
\begin{align}
    x_{i, k+1} & =x_{i, k}-\alpha\left[\nabla_1 f_i\left(x_{i, k}, u_{i, k}\right)+\nabla \phi_i\left(x_{i, k}\right) s_{i, k}\right] \notag\\
& \quad +\beta (x_{i,k}-x_{i,k-1}), \\
u_{i, k+1} & =\sum_{j=1}^N a_{i j} u_{j, k}+\phi_i\left(x_{i, k+1}\right)-\phi_i\left(x_{i, k}\right), \\
s_{i, k+1} & =\sum_{j=1}^N a_{i j} s_{j, k}+\nabla_2 f_i\left(x_{i, k+1}, u_{i, k+1}\right)  \notag\\
&\quad  -\nabla_2 f_i\left(x_{i, k}, u_{i, k}\right). 
\end{align}

In DAGT-HB, we introduce momentum term $x_{i,k}-x_{i,k-1}$ to accelerate the convergence rate of the algorithm. Because distributed aggregative optimization has departed from consistency protocols and coupled  agents together through an aggregative variable, each  agent only needs to reach the optimal point of its own local objective function. Hence, intuitively speaking, for each agent, if the iteration point generated by the algorithm always moves towards the optimal solution, adding a momentum in the same direction can inevitably accelerate the convergence speed of the algorithm. Because $u(x)$ is global information that cannot be accessed directly for all agents, $u_{i,k}$ is introduced for agent $i$ to track the average $u(x)$. Meanwhile, $s_{i,k}$ tracks the gradient sum $\frac{1}{N}\sum_{i=1}^N\nabla_2 f_i\left(x_{i}, u(x)\right)$ which cannot also be obtained to all agents.

The DAGT-HB algorithm can be rewritten as the following compact form:
\begin{align}
    x_{k+1} & =x_k-\alpha\left[\nabla_1 f\left(x_k, u_k\right)+\nabla \phi\left(x_k\right) s_k\right]+\beta(x_k-x_{k-1}), \\
u_{k+1} & =\mathcal{A} u_k+\phi\left(x_{k+1}\right)-\phi\left(x_k\right), \\
s_{k+1} & =\mathcal{A} s_k+\nabla_2 f\left(x_{k+1}, u_{k+1}\right)-\nabla_2 f\left(x_k, u_k\right),
\end{align}
with $\mathcal{A}=A \otimes I_d$ as defined in Lemma $2$, $x_k=$ col$\left(x_{1, k}, \ldots, x_{N, k}\right)$, and similar notations for $u_k$ and $s_k$.

Firstly we notice that multiplying $\frac{1}{N}$ on both sides of (11) and (12) can lead to
\begin{align}
    &\bar{u}_{k+1}=\bar{u}_k+\frac{1}{N} \sum_{i=1}^N \phi_i\left(x_{i, k+1}\right)-\frac{1}{N} \sum_{i=1}^N \phi_i\left(x_{i, k}\right), \\
&\bar{s}_{k+1}=\bar{s}_k+\frac{1}{N} \sum_{i=1}^N\left[\nabla_2 f_i\left(x_{i,k+1}, u_{i,k+1}\right) -\nabla_2 f_i\left(x_{i,k}, u_{i,k}\right)\right].
\end{align}
Then if  we initialize $u$ and $s$ as $u_{i,0}=\phi_i(x_{i,0})$ and $s_{i,0}=\nabla_2 f_i(x_{i,0},u_{i,0})$ for $i=1,2,\ldots,N$, where $x_{i,0}$ is arbitrary, we can derive 
\begin{align}
    \bar{u}_k & =\frac{1}{N} \sum_{i=1}^N u_{i, k}=\frac{1}{N} \sum_{i=1}^N \phi_i\left(x_{i, k}\right), \\
\bar{s}_k & =\frac{1}{N} \sum_{i=1}^N s_{i, k}=\frac{1}{N} \sum_{i=1}^N \nabla_2 f_i\left(x_{i, k}, u_{i, k}\right) .
\end{align}

Next, we establish the equivalence of the optimal solution to the problem (1) and the fixed point of the DAGT-HB algorithm. 

\textbf{Lemma 5:}  Under Assumption $1$ and Assumption $2$, the fixed point of (10)-(12) is
the optimal solution to problem (1).

\begin{proof}Denote the equilibrium point of (11)-(13) as $x^*=$col $\left(x_1^*, \ldots, x_N^*\right)$, $u^*=$ col$\left(u_1^*, \ldots, u_N^*\right)$, and $s^*=\operatorname{col}\left(s_1^*, \ldots, s_N^*\right)$. From (10)-(12), we can obtain
\begin{align}
    & \nabla_1 f\left(x^*, u^*\right)+\nabla \phi\left(x^*\right) s^*=\mathbf{0}_{Nd}, \\
& \mathcal{L} u^*=\mathbf{0}_{Nd},\quad \mathcal{L} s^*=\mathbf{0}_{Nd} ,
\end{align}
where $\mathcal{L}=L\otimes I_d$ and $L$ is Laplacian matrix of graph $\mathcal{G}$. Due to  the properties of the Laplace matrix , it is easy to derive that $u_i^*=u_j^*$ and $s_i^*=s_j^*$ for all $i\neq j$. Because of formulas (15)-(16), it leads to
\begin{align}
    u_i^* & =\frac{1}{N} \sum_{i=1}^N \phi_i\left(x_{i}^*\right)=u(x^*), \\
s_i^* & =\frac{1}{N} \sum_{i=1}^N \nabla_2 f_i\left(x_{i}^*, u(x^*)\right) .
\end{align}
By substituting (19) and (20) into (17), we can obtain $\nabla f(x^*)=0$. Because $f$ is $m-$ strongly convex, $x^*$ is the unique optimal solution to problem (1).
\end{proof}

\subsection{Auxiliary Results}
In order to analyze the convergence and convergence rate of this algorithm, we use the method of compressed state vector and collects the following four quantities:
\begin{enumerate}
    \item $||x_{k+1}-x^*||$, the state error in the network;
    \item $||x_{k+1}-x_{k}||$, the state difference;
    \item $||u_{k+1}-\mathcal{K}u_{k+1}||$, the aggregative variable tracking error;
    \item $||s_{k+1}-\mathcal{K}s_{k+1}||$,  the gradient sum  tracking error.
\end{enumerate}

In the next Lemmas 6–9, we derive the relationships among the four quantities mentioned above. Firstly, we derive the bound on  $||x_{k+1}-x^*||$, the state error in the network.

\textbf{Lemma 6:} Under Assumptions $1$-$4$, the  following inequality holds, $\forall k\geq  0$:
\begin{equation}
    \begin{aligned}
        & \left\|x_{k+1}-x^*\right\| \\
        & \leq(1-m \alpha)\left\|x_k-x^*\right\|+\alpha L_1\left\|u_k-\mathcal{K} u_k\right\|\\
&\quad+\alpha L_3\left\|s_k-\mathcal{K} s_k\right\|+\beta\left\|x_k-x_{k-1}\right\|.
    \end{aligned}
\end{equation}

\begin{proof}
For $\left\|x_{k+1}-x^*\right\|$, by invoking (10), it leads to
    \begin{equation}
        \begin{aligned}
            & \left\|x_{k+1}-x^*\right\| \\
& =\left\|x_k-x^*-\alpha\left[\nabla_1 f\left(x_k, u_k\right)+\nabla \phi\left(x_k\right) s_k\right]+\beta(x_k-x_{k-1})\right\| \\
& \leq \left\|x_k-x^*-\alpha\left[\nabla_1 f\left(x_k, u_k\right)+\nabla \phi\left(x_k\right) s_k\right]\right\|+\beta\left\|x_k-x_{k-1}\right\| \\
& \leq \Big\Vert x_k-x^*-\alpha\left[\nabla_1 f\left(x_k, \mathbf{1}_N \otimes \bar{u}_k\right)\right. \\
&\quad \left.+\nabla \phi\left(x_k\right) \mathbf{1}_N \otimes \frac{1}{N} \sum_{i=1}^N \nabla_2 f_i\left(x_{i, k}, \mathbf{1}_N \otimes \bar{u}_k\right)\right]+\alpha \nabla f\left(x^*\right) \Big\Vert \\
&\quad +\alpha \Big\Vert \nabla_1 f\left(x_k, u_k\right)+\nabla \phi\left(x_k\right) \mathbf{1}_N \otimes \bar{s}_k-\nabla_1 f\left(x_k, \mathbf{1}_N \otimes \bar{u}_k\right) \\
&\quad -\nabla \phi\left(x_k\right) \mathbf{1}_N \otimes \frac{1}{N} \sum_{i=1}^N \nabla_2 f_i\left(x_{i, k}, \mathbf{1}_N \otimes \bar{u}_k\right) \Big\Vert \\
&\quad +\alpha\left\|\nabla \phi\left(x_k\right) s_k-\nabla \phi\left(x_k\right) \mathbf{1}_N \otimes \bar{s}_k\right\|+\beta\left\|x_k-x_{k-1}\right\|. \\
        \end{aligned}
    \end{equation}
From Lemma $1$, we can bound the first term of the right term of (22) as follows:
\begin{equation}
    \begin{aligned}
        &\Big\Vert x_k-\alpha[\nabla \phi\left(x_k\right)
         \mathbf{1}_N \otimes \frac{1}{N} \sum_{i=1}^N \nabla_2 f_i\left(x_{i, k}, \mathbf{1}_N \otimes \bar{u}_k\right) \\
         &\quad+\nabla_1 f\left(x_k, \mathbf{1}_N \otimes \bar{u}_k\right)]-\left[x^*-\alpha \nabla f\left(x^*\right)\right]\Big\Vert \\
         &  \leq(1-m \alpha)\left\|x_k-x^*\right\|.
    \end{aligned}
\end{equation}
For the second term, since $f(x)$ is $L_1-$smooth and $\mathbf{1}_N \otimes \bar{u}_k=\mathcal{K} u_k$ we can obtain
\begin{equation}
    \begin{aligned}
        & \alpha \Big\Vert \nabla_1 f\left(x_k, u_k\right)+\nabla \phi\left(x_k\right) \mathbf{1}_N \otimes \bar{s}_k-\nabla_1 f\left(x_k, \mathbf{1}_N \otimes \bar{u}_k\right) \\
& \quad-\nabla \phi\left(x_k\right) \mathbf{1}_N \otimes \frac{1}{N} \sum_{i=1}^N \nabla_2 f_i\left(x_{i, k}, \mathbf{1}_N \otimes \bar{u}_k\right) \Big\Vert \\
& \leq \alpha L_1\left\|u_k-\mathcal{K} u_k\right\|.
    \end{aligned}
\end{equation}
For the third term, using Assumption $4$ and $\mathbf{1}_N \otimes \bar{s}_k=\mathcal{K} s_k$ can get the following inequality:
\begin{equation}
    \begin{aligned}
        \alpha\left\|\nabla \phi\left(x_k\right) s_k-\nabla \phi\left(x_k\right) \mathbf{1}_N \otimes \bar{s}_k\right\| \leq \alpha L_3\left\|s_k-\mathcal{K} s_k\right\|.
    \end{aligned}
\end{equation}
Then by substituting (23)-(25) into (22), we complete the proof.
\end{proof}

Secondly, we derive a bound for $\left\|x_{k+1}-x_k\right\|$.

\textbf{Lemma $7$:} Under Assumptions $1$-$4$, the  following inequality holds, $\forall k\geq  0$:
\begin{equation}
    \begin{aligned}
\left\|x_{k+1}-x_k\right\| \leq & \alpha L_1\left(1+L_3\right)\left\|x_k-x^*\right\|+\alpha L_1\left\|u_k-\mathcal{K} u_k\right\| \\
& +\alpha L_3\left\|s_k-\mathcal{K} s_k\right\|+\beta\left\|x_k-x_{k-1}\right\|.
\end{aligned}
\end{equation}

\begin{proof}
Note that $\nabla f(x^*)=0$ and then we have
    \begin{equation}
        \begin{aligned}
            & \left\|x_{k+1}-x_k\right\| \\
& =\left\|\beta(x_k-x_{k-1})-\alpha(\nabla_1 f\left(x_k, u_k\right)+\nabla \phi\left(x_k\right) s_k)\right\| \\
& \leq \alpha \Bigg\Vert \nabla_1 f\left(x_k, u_k\right)+\nabla \phi\left(x_k\right) \mathcal{K} s_k-\nabla_1 f\left(x^*, \mathbf{1}_N \otimes u^*\right) \\
& \quad\quad-\nabla \phi\left(x^*\right)\left[\mathbf{1}_N \otimes \frac{1}{N} \sum_{i=1}^N \nabla_2 f_i\left(x^*, \mathbf{1}_N \otimes u^*\right)\right] \Bigg\Vert \\
& \quad+\alpha\left\|\nabla \phi\left(x_k\right)\left(s_k-\mathcal{K} s_k\right)\right\|+\beta\left\|x_k-x_{k-1}\right\|. \\
        \end{aligned}
    \end{equation}
By utilizing Assumption $2$ and triangle inequality of norm, we can obtain the following formula:
\begin{equation}
    \begin{aligned}
       &   \Big\Vert \nabla_1 f\left(x_k, u_k\right)+\nabla \phi\left(x_k\right) \mathcal{K} s_k-\nabla_1 f\left(x^*, \mathbf{1}_N \otimes u^*\right) \\
& \quad-\nabla \phi\left(x^*\right)\left[\mathbf{1}_N \otimes \frac{1}{N} \sum_{i=1}^N \nabla_2 f_i\left(x^*, \mathbf{1}_N \otimes u^*\right)\right] \Big\Vert \\ 
& \leq L_1\left(\left\|x_k-x^*\right\|+\left\|u_k-\mathbf{1}_N \otimes u^*\right\|\right)\\
& \leq L_1\left(\left\|x_k-x^*\right\|+\left\|u_k-\mathcal{K}u_k\right\|+\left\|\mathcal{K}u_k-\mathbf{1}_N \otimes u^*\right\|\right).
    \end{aligned}
\end{equation}
For $\left\|\mathcal{K}u_k-\mathbf{1}_N \otimes u^*\right\|$, we can derive
\begin{equation}
    \begin{aligned}
\left\|\mathcal{K} u_k-1_N \otimes u^*\right\|^2 & =\left\|1_N \otimes\left(\bar{u}_k-u^*\right)\right\|^2 \\
& =N\left\|\frac{1}{N} \sum_{i=1}^N\left(\phi_i\left(x_{i, k}\right)-\phi_i\left(x_i^*\right)\right)\right\|^2 \\
& \leq \frac{1}{N}\left(\sum_{i=1}^N\left\|\phi_i\left(x_{i, k}\right)-\phi_i\left(x_i^*\right)\right\|\right)^2 \\
& \leq \frac{1}{N}\left(\sum_{i=1}^N L_3\left\|x_{i, k}-x_i^*\right\|\right)^2 \\
& \leq L_3^2 \sum_{i=1}^N\left\|x_{i, k}-x_i^*\right\|^2 \\
& =L_3^2\left\|x_k-x^*\right\|^2,
\end{aligned}
\end{equation}
where using the property that $\phi_i$ is $L_3$-Lipschitz continuous can obtain the second inequality, and applying the fact that $\left(\sum_{i=1}^N a_i\right)^2 \leq$ $N \sum_{i=1}^N a_i^2$ for any nonnegative scalars $a_i$s  can easily get the last inequality. Thus, we have
\begin{equation}
    \left\|\mathcal{K} u_k-1_N \otimes u^*\right\| \leq L_3\left\|x_k-x^*\right\|.
\end{equation}
Then by using Assumption $4$ we can obtain
\begin{equation}
    \left\|\nabla \phi\left(x_k\right)\left(s_k-\mathcal{K} s_k\right)\right\| \leq L_3 \left\|s_k-\mathcal{K} s_k\right\|.
\end{equation}
Finally, by inserting (28), (30) and (31) into (27) we can obtain the result (26).
\end{proof}

The next step is to bound  the aggregative variable tracking error $||u_{k+1}-\mathcal{K}u_{k+1}||$.

\textbf{Lemma $8$:} Under Assumptions $1$-$4$, the  following inequality holds, $\forall k\geq  0$:
\begin{equation}
\begin{aligned}
&\left\|u_{k+1}-\mathcal{K} u_{k+1}\right\| \\
& \leq\left(\rho+\alpha L_1 L_3\right)\left\|u_k-\mathcal{K} u_k\right\|+\alpha L_1 L_3\left(1+L_3\right)\left\|x_k-x^*\right\| \\
&\quad+\alpha L_3^2\left\|s_k-\mathcal{K} s_k\right\|+\beta L_3\left\|x_k-x_{k-1}\right\|.
\end{aligned}
\end{equation}

\begin{proof}For $\left\|u_{k+1}-\mathcal{K} u_{k+1}\right\|$, by invoking (11), it leads to
    \begin{equation}
\begin{aligned}
& \left\|u_{k+1}-\mathcal{K} u_{k+1}\right\| \\
& =\left\|\mathcal{A} u_k+\phi\left(x_{k+1}\right)-\phi\left(x_k\right)-\mathcal{K} \mathcal{A} u_k-\mathcal{K}\left[\phi\left(x_{k+1}\right)-\phi\left(x_k\right)\right]\right\| \\
& \leq \rho\left\|u_k-\mathcal{K} u_k\right\|+\|I-\mathcal{K}\|\left\|\phi\left(x_{k+1}\right)-\phi\left(x_k\right)\right\| \\
& \leq \rho\left\|u_k-\mathcal{K} u_k\right\|+L_3\|I-\mathcal{K}\|\left\|x_{k+1}-x_k\right\|,
\end{aligned}
\end{equation}
where Lemma $2$ has been utilized to obtain the first inequality, and by using  Assumption $4$ we can obtain  the last inequality.  Notice that $\|I-\mathcal{K}\|=1$, then by substituting (26) into (33) we can complete the proof.
\end{proof}

Lastly, we derive the bound  $||s_{k+1}-\mathcal{K}s_{k+1}||$,  the gradient sum $\frac{1}{N}\sum_{i=1}^N\nabla_2 f_i\left(x_{i}, u(x_{k+1})\right)$ tracking error.

\textbf{Lemma $9$:} Under Assumptions $1$-$4$, the  following inequality holds, $\forall k\geq  0$:
\begin{equation}
    \begin{aligned}
        & \left\|s_{k+1}-\mathcal{K} s_{k+1}\right\| \\
& \leq \left(\rho+\alpha L_2 L_3\left(1+L_3\right)\right)\left\|s_k-\mathcal{K} s_k\right\| 
  \\
& \quad+\left(\alpha L_1 L_2\left(1+L_3\right)+2L_2\right)\left\|u_k-\mathcal{K} u_k\right\|  \\
&\quad+ \alpha L_1 L_2\left(1+L_3\right)^2\left\|x_k-x^*\right\|+\beta L_2(1+L_3)||x_k-x_{k-1}||.
    \end{aligned}
\end{equation}

\begin{proof}For $\left\|s_{k+1}-\mathcal{K} s_{k+1}\right\|$, by invoking (13), it leads to
    \begin{equation}
\begin{aligned}
& \left\|s_{k+1}-\mathcal{K} s_{k+1}\right\| \\
& =||\mathcal{A} s_k+\nabla_2 f\left(x_{k+1}, u_{k+1}\right)-\nabla_2 f\left(x_k, u_k\right)\\
&\quad\quad-\mathcal{K}\mathcal{A} s_k-\mathcal{K}[\nabla_2 f\left(x_{k+1}, u_{k+1}\right)-\nabla_2 f\left(x_k, u_k\right)]|| \\
& \leq \left\|\mathcal{A}s_k-\mathcal{K} s_k\right\|\\
& \quad+\|I-\mathcal{K}\|\left\|\nabla_2 f\left(x_{k+1}, u_{k+1}\right)-\nabla_2 f\left(x_k, u_k\right)\right\|\\
& \leq \rho\left\|s_k-\mathcal{K} s_k\right\|+\left\|\nabla_2 f\left(x_{k+1}, u_{k+1}\right)-\nabla_2 f\left(x_k, u_k\right)\right\| \\
& \leq \rho\left\|s_k-\mathcal{K} s_k\right\|+  L_2\left(\left\|x_{k+1}-x_k\right\|+\left\|u_{k+1}-u_k\right\|\right),\\
\end{aligned}
\end{equation}
where  Lemma $2$ has been utilized to obtain the first inequality and Assumption $2$ has been leveraged in the last inequality.  Notice that
\begin{equation}
\begin{aligned}
    & \left\|u_{k+1}-u_k\right\|\\
    & =\left\|\mathcal{A}u_k-u_k+\phi\left(x_{k+1}\right)-\phi\left(x_k\right)\right\|\\
    & =\left\|\mathcal{A}u_k-\mathcal{A}\mathcal{K}u_k+\mathcal{K}u_k-u_k+\phi\left(x_{k+1}\right)-\phi\left(x_k\right)\right\|\\
    & \leq \|(\mathcal{A}-\left.I_N \otimes I_d\right)\left(u_k-\mathcal{K} u_k\right)\|+\|\phi\left(x_{k+1}\right)-\phi\left(x_k\right)\|\\
    & \leq ||A-I_N|| \left\|u_k-\mathcal{K} u_k\right\|+ L_3\left\|x_{k+1}-x_k\right\| \\
    & \leq 2 \left\|u_k-\mathcal{K} u_k\right\|+ L_3\left\|x_{k+1}-x_k\right\|. \\
    \end{aligned}
\end{equation}
Then by substituting  (36) and (26) into (35), we can finish the proof.
\end{proof}

\subsection{Main Result}
We now present the main result of this section. Based on Lemmas $6$-$9$, we give the convergence and convergence rate of the DAGT-HB algorithm in the following theorem.

\textbf{Theorem $1$}: Under Assumptions $1-4$, if 
\begin{equation}
    (\alpha,\beta) \in \bigcap\limits _{i=1}^6\mathcal{S}_i,
\end{equation}
where $\mathcal{S}_i,i=1,\cdots,6$  are defined in the following proof, then $x_k=$ col$\left(x_{1, k}, \ldots, x_{N, k}\right)$ generated by DAGT-HB can converge to the
optimizer of problem (1) at the $R-$linear convergence rate.

\begin{proof}
 Denote
\begin{equation}
    V_k=\text{col}\left(\left\|x_k-x^*\right\|,\left\|x_k-x_{k-1}\right\|,\left\|u_k-\mathcal{K} u_k\right\|,\left\|s_k-\mathcal{K} s_k\right\|\right).
\end{equation}
From Lemmas $6$-$9$, it can be concluded that
\begin{equation}
    V_{k+1} \leq P V_k,
\end{equation}
where 
\begin{equation}
     \begin{aligned}
     P=\left[\begin{array}{cccc}
1-m\alpha & \beta & \alpha L_1 & \alpha L_3\\
\alpha L_1(1+L_3) & \beta &\alpha L_1 & \alpha L_3\\
\alpha L_1L_3(1+L_3) & \beta L_3 & \rho+\alpha L_1L_3 & \alpha L_3^2\\
\alpha L_1L_3(1+L_3)^2 & \beta L_2(1+L_3) & p_{43} & p_{44}
\end{array}\right],\\
     \end{aligned}
 \end{equation}
 $p_{43}=\alpha L_1L_3(1+L_3)+2L_2$ and $p_{44}=\rho+\alpha L_2L_3(1+L_3).$
Firstly, based on Lemma $3$, we seek a small range of $\alpha$ and $\beta$
 to satisfy $\rho(P)<1$. We define a positive vector $z=[z_1,z_2,z_3,z_4]^T$ such that
 \begin{equation}
     Pz<z,
 \end{equation}
 which is equal to
 \begin{equation}
     \begin{aligned}
        & 0<\beta<\frac{\alpha(m z_1-L_1z_3-L_3z_4)}{z_2}=M_1,\\
       &  0<\beta<\frac{z_2-\alpha L_1(1+L_3)z_1-\alpha L_1z_3-\alpha L_3 z_4}{z_2}=M_2,\\
       &  0<\beta<\\
       &\frac{(1-\rho-\alpha L_1L_3)z_3-\alpha L_1L_3(1+L_3)z_1-\alpha L_3^2z_4}{L_3z_2}=M_3,\\
       & 0<\beta<\frac{[1-\rho-\alpha L_2L_3(1+L_3)]z_4-\alpha L_1L_2(1+L_3)^2z_1}{L_2(1+L_3)z_2},\\
       &- \frac{[\alpha L_1L_2(1+L_3)+2L_2]z_3}{L_2(1+L_3)z_2}=M_4.
     \end{aligned}
 \end{equation}
 From the above inequalities, we derive
 \begin{equation}
     \begin{aligned}
         & 0< \alpha <\frac{z_2}{L_1(1+L_3)z_1+L_1z_3+L_3 z_4}=J_1,\\
         & 0< \alpha <\frac{1-\rho}{L_1L_3}=J_2,\\
         & 0< \alpha <\frac{(1-\rho)z_3}{L_1L_3(1+L_3)z_1+L_1L_3z_3+L_3^2z_4},\\
         & 0< \alpha < \frac{1-\rho}{ L_2L_3(1+L_3)}=J_3,\\
         & 0< \alpha < 
         \frac{(1-\rho)z_4-2L_2z_3}{L_2(1+L_3)[L_1(1+L_3)z_1+L_1z_3+L_3z_4]}=J_4,\\
         & z_1>\frac{L_1z_3+L_3z_4}{m},\\
         & z_4>\frac{2L_2z_3}{1-\rho}.
     \end{aligned}
 \end{equation}
 That is to say, we can select arbitrary $z_2$ and $z_3$, when
 \begin{equation}
     \begin{aligned}
       &  0<\alpha<\bar{\alpha}=\min \{ J_1, J_2, J_3, J_4, \frac{1}{L_1} \},
     \end{aligned}
 \end{equation}
 and 
 \begin{equation}
     \begin{aligned}
         & 0<\beta<\bar{\beta}=\min \{M_1, M_2, M_3, M_4\},
     \end{aligned}
 \end{equation}
 where $z_2>0, z_3>0, z_1>\frac{L_1z_3+L_3z_4}{\mu}$ and $z_4>\frac{2L_2z_3}{1-\rho}$, we have $\rho(P)<1$. 

 Next, we use Jury criterion to seek precise range of $\alpha$ and $\beta$ to meet $\rho(P)<1$. By computing, we can obtain the characteristic polynomial of $P$:
 \begin{equation}
     H(\lambda)=|\lambda I-P|=a_0+a_1\lambda+a_2\lambda^2+a_3\lambda^3+\lambda^4,
 \end{equation}
 where
 \begin{equation}
     \begin{aligned}
         a_0 &=\beta d_1 \rho(\rho+2 \alpha d_2),\\
a_1&=  \beta[-d_1 \rho+(d_1+\rho)\times
\left(-\rho+2 \alpha d_2\right)] 
 +(m \alpha-1) \rho\left(\rho+\alpha d_2\right)\\
&\quad-\alpha L_3 d_1 \times
[\alpha L_1\left(\rho+
 \alpha L_3 d_2\right)+\alpha L_3 d_3 ],\\
 a_2&=\beta(d_1+2\rho+\alpha d_2)+(1-m\alpha)(2\rho+\alpha d_2)+\rho(\rho+\alpha d_2)\\
 &\quad+\alpha L_3 d_1[L_1+L_3(1+L_3)]+L_3[\alpha L_1(\rho+\alpha d_2)+\alpha L_3 d_3],\\
 a_3&=-\beta+(m-2)\alpha-1-\alpha L_3[L_2(1+L_3)+L_1],
     \end{aligned}
 \end{equation}
 and $d_1=1-m\alpha-\alpha L_1(1+L_3)$, $d_2=L_3(1+L_3)(L_2-L_3)$ and $d_3=-2\alpha L_2+\rho(1+L_3)$. Then we can obtain:
 \begin{equation}
     \begin{aligned}
         &b_0=a_0^2-1, b_1=a_0a_1-a_3, b_2=a_0a_2-a_2, b_3=a_0a_3-a_1\\
         &c_0=b_0^2-b_3^2,c_1=b_0b_1-b_2b_3,c_2=b_0b_2-b_1b_3.
     \end{aligned}
 \end{equation}
 Next we denote
 \begin{equation}
     \begin{aligned}
         &\mathcal{S}_1=\{(\alpha,\beta)|a_0+a_1+a_2+a_3+1>0\},\\
        &\mathcal{S}_2=\{(\alpha,\beta)|a_0-a_1+a_2-a_3+1>0\},\\
        &\mathcal{S}_3=\{(\alpha,\beta)||a_0|<1\},\quad\mathcal{S}_4=\{(\alpha,\beta)||b_0|>|b_3|\}\\
        &\mathcal{S}_5=\{(\alpha,\beta)||c_0|>|c_2|\},\quad\mathcal{S}_6=\{(\alpha,\beta)|\alpha>0,\beta>0\}.
     \end{aligned}
 \end{equation}
 Hence according to Lemma $4$, when $(\alpha,\beta) \in \bigcap\limits _{i=1}^6\mathcal{S}_i$, the spectral radius of the matrix P is less than 1. In view of the above analysis, we know $\bigcap\limits _{i=1}^6\mathcal{S}_i$ is non-empty.
 Finally, denote $\rho_1=\rho(P)$ and $0<\rho_1<1$ then we can obtain 
 \begin{equation}
     ||V_{k+1}||\leq \rho_1||V_k||.
 \end{equation}
 Furthermore, it leads to
 \begin{equation}
     ||x_k-x^*||\leq ||V_k|| \leq C_1\rho_1^k,
 \end{equation}
 where $C_1=||V_0||$. So DAGT-HB can achieve the
$\mathbf{R}$-linear convergence rate. Then the proof is completed.
\end{proof}

\textbf{Remark $1$:} In Theorem $1$, we have established an $\mathbf{R}$-linear rate of DAGT-HB when the step-size ${\alpha}$, and the  momentum
parameter ${\beta}$ follow (37). But we acknowledge
that the theoretical bounds of ${\alpha}$ and ${\beta}$ in Theorem 1 are conservative. How to obtain theoretical boundaries and even optimal parameters will be considered in our future work.

\textbf{Corollary $1$:} Under the same assumptions of
Theorem $1$, the following equality holds:
\begin{equation}
    |f(x_k)-f(x^*)|\leq \frac{L_1}{2}C_1^2 \rho_1^{2k}.
\end{equation}
\begin{proof}Because $f(x)$ is $L_1$-smooth and $\nabla f(x^*)=0$, we can obtain
\begin{equation}
    f(x_k)-f(x^*)\leq \frac{L_1}{2}||x_k-x^*||^2.
\end{equation}
By substituting (51) into (53), we complete the proof.
\end{proof}

As a summary of this section, the DAGT-HB method is formulated as the following Algorithm \ref{algorithm1}.

\begin{algorithm}
    \caption{DAGT-HB}\label{algorithm1}
    \SetKwInOut{Input}{Input}\SetKwInOut{Output}{Output}
    \Input{initial point $x_{i,-1} $ and $x_{i,0} \in \mathbb{R}^{n_i}$, $u_{i,0}=\phi_i(x_{i,0})$ and $s_{i,0}=\nabla_2 f_i(x_{i,0},u_{i,0})$ for $i=1,2,\ldots,N$, let $  (\alpha,\beta) \in \bigcap\limits _{i=1}^6\mathcal{S}_i$, set $k=0$, select appropriate $\epsilon>0$ and $K_{max} \in \mathbb{N}^+$. }
    \Output{optimal $x^*$ and $f(x^*)$.}
    While
    {$k<K_{max}$ and $||\nabla f(x_k)||\geq \epsilon$ do}\\
    Iterate: Update for each $i \in \{1, 2, \ldots, N\}$:
    \begin{align*}
        x_{i,k+1}\;\text{iterated by formula (7),}\\
        u_{i,k+1}\;\text{iterated by formula (8),}\\
        s_{i,k+1}\;\text{iterated by formula (9).}
    \end{align*}
    
    Update: $k=k+1$.\\
\end{algorithm}

\section{DAGT-NES}
In addition to heavy ball, the Nesterov’s algorithm is also a well-known accelerated  method that can be combined with DAGT to solve problem (1). To this end, we propose the following DAGT-NES algorithm:
\begin{align}
  &  x_{i, k+1}  =y_{i, k}-\alpha\left[\nabla_1 f_i\left(y_{i, k}, u_{i, k}\right)+\nabla \phi_i\left(y_{i, k}\right) s_{i, k}\right], \\
 &y_{i,k+1} = x_{i,k+1} + \gamma(x_{i,k+1}-x_{i,k}),\\
&u_{i, k+1}  =\sum_{j=1}^N a_{i j} u_{j, k}+\phi_i\left(y_{i, k+1}\right)-\phi_i\left(y_{i, k}\right), \\
&s_{i, k+1} =\sum_{j=1}^N a_{i j} s_{j, k}+\nabla_2 f_i\left(y_{i, k+1}, u_{i, k+1}\right) -\nabla_2 f_i\left(y_{i, k}, u_{i, k}\right). 
\end{align}
Likewise, in DAGT-NES, $u_{i,k}$ is introduced for agent $i$ to track the average $u(y)$ and $s_{i,k}$ tracks the gradient sum $\frac{1}{N}\sum_{i=1}^N\nabla_2 f_i\left(y_{i}, u(y)\right)$ .
 The DAGT-NES algorithm can be rewritten as the following compact form:
\begin{align}
    x_{k+1} & =y_k-\alpha\left[\nabla_1 f\left(y_k, u_k\right)+\nabla \phi\left(y_k\right) s_k\right], \\
    y_{k+1} & = x_{k+1}+\gamma (x_{k+1}-x_k),\\
u_{k+1} & =\mathcal{A} u_k+\phi\left(y_{k+1}\right)-\phi\left(y_k\right), \\
s_{k+1} & =\mathcal{A} s_k+\nabla_2 f\left(y_{k+1}, u_{k+1}\right)-\nabla_2 f\left(y_k, u_k\right).
\end{align}
with $\mathcal{A}=A \otimes I_d$ as defined in Lemma $2$, $x_k=$ col$\left(x_{1, k}, \ldots, x_{N, k}\right)$, and similar notations for $u_k, y_k$ and $s_k$.

Note that if  we initialize $u$ and $s$ as $u_{i,0}=\phi_i(y_{i,0})$ and $s_{i,0}=\nabla_2 f_i(y_{i,0},u_{i,0})$ for $i=1,2,\ldots,N$, where $y_{i,0}$ is arbitrary, analogous to DAGT-HB, we can obtain
\begin{align}
    \bar{u}_k & =\frac{1}{N} \sum_{i=1}^N u_{i, k}=\frac{1}{N} \sum_{i=1}^N \phi_i\left(y_{i, k}\right), \\
\bar{s}_k & =\frac{1}{N} \sum_{i=1}^N s_{i, k}=\frac{1}{N} \sum_{i=1}^N \nabla_2 f_i\left(y_{i, k}, u_{i, k}\right) .
\end{align}

Next, we establish the equivalence of the optimal solution to the problem (1) and the fixed point of the DAGT-NES algorithm. 

\textbf{Lemma 10}: Under Assumption $1$ and Assumption $2$, the equilibrium point of (58)-(61) is the optimal solution to problem (1).

\begin{proof}See Appendix A.
\end{proof}

\subsection{Auxiliary Results}
Similar to DAGT-HB, we utilize  the method of compressed state vector to derive the convergence and convergence rate of DAGT-NES and still collect the following four quantities:
\begin{enumerate}
    \item $||x_{k+1}-x^*||$, the state error in the network;
    \item $||x_{k+1}-x_{k}||$, the state difference;
    \item $||u_{k+1}-\mathcal{K}u_{k+1}||$, the aggregative variable tracking error;
    \item $||s_{k+1}-\mathcal{K}s_{k+1}||$,  the gradient sum  tracking error;
\end{enumerate}

In the next Lemmas 11–14, we derive the relationship among the four quantities mentioned above. Firstly, we derive the bound on  $||x_{k+1}-x^*||$, the state error in the network.

\textbf{Lemma 11:} Under Assumptions $1$-$4$, the  following inequality holds, $\forall k\geq  0$:
\begin{equation}
    \begin{aligned}
        & \left\|x_{k+1}-x^*\right\| \\
        & \leq(1-m \alpha)\left\|x_k-x^*\right\|+\alpha L_1\left\|u_k-\mathcal{K} u_k\right\|\\
&\quad+\alpha L_3\left\|s_k-\mathcal{K} s_k\right\|+(1-m \alpha)\gamma\left\|x_k-x_{k-1}\right\|.\\
    \end{aligned}
\end{equation}

\begin{proof}See Appendix B.
\end{proof}

Secondly, we derive a bound for $\left\|x_{k+1}-x_k\right\|$.

\textbf{Lemma 12:} Under Assumptions $1-4$, the  following inequality holds, $\forall k\geq  0$:
\begin{equation}
    \begin{aligned}
        &\left\|x_{k+1}-x_k\right\|\\
        &\leq  \alpha L_1\left(1+L_3\right)\left\|x_k-x^*\right\|+\alpha L_1\left\|u_k-\mathcal{K} u_k\right\| \\
& \quad+\alpha L_3\left\|s_k-\mathcal{K} s_k\right\|+\gamma(1+\alpha L_1+\alpha L_1L_3)\left\|x_k-x_{k-1}\right\|.
    \end{aligned}
\end{equation}
\begin{proof}See Appendix C.
\end{proof}

The next step is to bound  the aggregative variable tracking error $||u_{k+1}-\mathcal{K}u_{k+1}||$.

\textbf{Lemma 13:} Under Assumptions $1$-$4$, the  following inequality holds, $\forall k\geq  0$:
\begin{equation}
    \begin{aligned}
       &\left\|u_{k+1}-\mathcal{K} u_{k+1}\right\| \\
& \leq\left[\rho+\alpha L_1 L_3(\gamma+1)\right]\left\|u_k-\mathcal{K} u_k\right\|+\alpha L_3^2(\gamma+1)\left\|s_k-\mathcal{K} s_k\right\|\\
&\quad +\alpha L_1 L_3\left(1+L_3\right)(\gamma+1)\left\|x_k-x^*\right\| \\
&\quad +\gamma L_3[(1+\gamma)(1+\alpha L_1 +\alpha L_1L_3)+1]\left\|x_k-x_{k-1}\right\|.\\ 
    \end{aligned}
\end{equation}

\begin{proof}See Appendix D.
\end{proof}

Lastly, we derive the bound  $||s_{k+1}-\mathcal{K}s_{k+1}||$,  the gradient sum $\frac{1}{N}\sum_{i=1}^N\nabla_2 f_i\left(x_{i}, u(x_{k+1})\right)$ tracking error.

\textbf{Lemma 14:} Under Assumptions $1$-$4$, the  following inequality holds, $\forall k\geq  0$:
\begin{equation}
    \begin{aligned}
        & \left\|s_{k+1}-\mathcal{K} s_{k+1}\right\| \\
& < \left[\rho+\alpha L_2 L_3\left(1+L_3\right)(1+\gamma)\right]\left\|s_k-\mathcal{K} s_k\right\| \\
&\quad +\alpha L_1 L_2\left(1+L_3\right)^2(\gamma+1)\left\|x_k-x^*\right\| \\
& \quad+\left[\alpha L_1 L_2\left(1+L_3\right)(r+1)+2L_2\right]\left\|u_k-\mathcal{K} u_k\right\|  \\
&\quad +\gamma L_2(L_3+1)[(1+\gamma)(1+\alpha L_1 +\alpha L_1L_3)+1]||x_k-x_{k-1}||.
    \end{aligned}
\end{equation}
\begin{proof}See Appendix E.
\end{proof}

\subsection{Main Result}
Summarizing 
Lemmas $11$-$14$, we give the convergence and convergence rate  of the DAGT-NES algorithm in the following theorem.

\textbf{Theorem $2$:} Under Assumptions 1-4, if 
\begin{equation}
   (\alpha,\gamma) \in \bigcap\limits _{i=1}^6\mathcal{T}_i,
\end{equation}
where $\mathcal{T}_i$, $i=1,\cdots,6$  are defined in the following proof, then $x_k=$ col$\left(x_{1, k}, \ldots, x_{N, k}\right)$ generated by DAGT-NES can converge to the
optimizer of problem (1) at the $R-$linear convergence rate.

\begin{proof}
    Denote
\begin{equation}
    V_k=\text{col}\left(\left\|x_k-x^*\right\|,\left\|x_k-x_{k-1}\right\|,\left\|u_k-\mathcal{K} u_k\right\|,\left\|s_k-\mathcal{K} s_k\right\|\right).
\end{equation}
From Lemmas 11-14, it can be concluded that
\begin{equation}
    V_{k+1} \leq Q V_k,
\end{equation}
where 
\begin{equation}
     \begin{aligned}
     Q=\left[\begin{array}{cccc}
1-m\alpha & (1-m\alpha)\gamma & \alpha L_1 & \alpha L_3\\
\alpha L_1(1+L_3) & q_{22} &\alpha L_1 & \alpha L_3\\
\alpha L_1 L_3\left(1+L_3\right)(\gamma+1) & q_{32} & q_{33} & q_{34} \\
\alpha L_1L_2(1+L_3)^2 (\gamma+1)& q_{42} & q_{43} & q_{44}
\end{array}\right],\\
     \end{aligned}
 \end{equation}
and 
\begin{equation}
\left\{\begin{array}{l}
 q_{22}=\gamma(1+\alpha L_1+\alpha L_1L_3),\\
 q_{32}=\gamma L_3[(1+\gamma)(1+\alpha L_1 +\alpha L_1L_3)+1],\\
q_{33}=\rho+\alpha L_1 L_3(\gamma+1),\\
q_{34}=\alpha L_3^2(\gamma+1),\\
q_{42}=\gamma L_2(L_3+1)[(1+\gamma)(1+\alpha L_1 +\alpha L_1L_3)+1],\\
q_{43}=\alpha L_1 L_2\left(1+L_3\right)(\gamma+1)+2L_2,\\
q_{44}=\rho+\alpha L_2 L_3\left(1+L_3\right)(1+\gamma).
\end{array}\right.
\end{equation}
Firstly, based on Lemma $3$, we want to seek a small range of $\alpha$ and $\gamma$
 to satisfy $\rho(Q)<1$. We define a positive vector $t=[t_1,t_2,t_3,t_4]^T$ such that
 \begin{equation}
     Qt<t .
 \end{equation}
However, since in the matrix $Q$, $\alpha$ and $\gamma$ have strong nonlinear relationship,  it is difficult to give the  range of $\alpha$ and $\gamma$ that makes the spectral radius of the matrix $Q$ less than 1. Thus,  to simplify the calculation, we let 
\begin{equation}
    \alpha\leq \frac{1}{L_1} \quad \text{and} \quad \gamma\leq \min \{\frac{1}{L_2},\frac{1}{L_3}\}
\end{equation}
 to eliminate some entries in matrix $Q$ that contain the nonlinear relationship between $\alpha$ and $\gamma$. Then we can obtain 
 \begin{equation}
    V_{k+1}<R V_k,
\end{equation}
where 
\begin{equation}
     \begin{aligned}
     R=\left[\begin{array}{cccc}
1-m\alpha & (1-m\alpha)\gamma & \alpha L_1 & \alpha L_3\\
\alpha L_1(1+L_3) & \gamma(2+L_3) &\alpha L_1 & \alpha L_3\\
\alpha L_1 L_3\left(2+L_3\right) & r_{32} & r_{33} & \alpha L_3(L_3+1)\\
r_{41} & r_{42} & r_{43} & r_{44}
\end{array}\right],\\
     \end{aligned}
 \end{equation}
and
\begin{equation}
\left\{\begin{array}{l}
 r_{32}=\gamma(L_3^2+4L_3+2),\\
r_{33}=\rho+\alpha L_1(L_3+1),\\
r_{41}=\alpha L_1(1+L_2)(1+L_3)^2,\\
r_{42}=\gamma(L_3+1)(L_2L_3+2L_2+L_3+1),\\
r_{43}=\alpha L_1 (L_2+1)\left(1+L_3\right)+2L_2,\\
r_{44}=\rho+\alpha L_2 \left(1+L_3\right)^2.
\end{array}\right.
\end{equation}
Then we solve the following equality:
\begin{equation}
    Rt<t,
\end{equation}
which is equivalent to
\begin{equation}
     \begin{aligned}
        & 0<\gamma<\frac{\alpha(m t_1-L_1t_3-L_3t_4)}{(1-m\alpha)t_2}=\Gamma_1,\\
       &  0<\gamma<\frac{t_2-\alpha L_1(1+L_3)t_1-\alpha L_1t_3-\alpha L_3 t_4}{(2+L_3)t_2}=\Gamma_2,\\
       &  0<\gamma<\frac{[1-\rho-\alpha L_1(L_3+1)]t_3-\alpha L_1L_3(2+L_3)t_1}{(L_3^2+4L_3+2)t_2}\\
       & -\frac{\alpha L_3(L_3+1)t_4}{(L_3^2+4L_3+2)t_2}=\Gamma_3,\\
       & 0<\gamma<\frac{[1-\rho-\alpha L_2(1+L_3)^2]t_4-\alpha L_1(L_2+1)(1+L_3)^2t_1}{L_2(1+L_3)t_2}\\
       &- \frac{[\alpha L_1(L_2+1)(1+L_3)+2L_2]t_3}{L_2(1+L_3)t_2}=\Gamma_4.
     \end{aligned}
 \end{equation}
 From the above inequalities, we derive
 \begin{equation}
     \begin{aligned}
         & 0< \alpha <\frac{t_2}{L_1(1+L_3)t_1+L_1t_3+L_3 t_4}=\Theta_1,\\
         & 0< \alpha <\frac{1-\rho}{L_1(L_3+1)}=\Theta_2,\\
         & 0< \alpha <\frac{(1-\rho)t_3}{L_1L_3(2+L_3)t_1+L_1(L_3+1)t_3+(L_3^2+L_3)t_4}=\Theta_3,\\
         & 0< \alpha < \frac{1-\rho}{ L_2L_3(1+L_3)}=\Theta_4,\\
         & 0< \alpha < \\
         &\frac{(1-\rho)t_4-2L_2t_3}{L_1(L_2+1)(1+L_3)[(1+L_3)t_1+t_3]+L_2(1+L_3)^2t_4}=\Theta_5,\\
         & t_1>\frac{L_1t_3+L_3t_4}{m},\\
         & t_4>\frac{2L_2t_3}{1-\rho}.
     \end{aligned}
 \end{equation}
 That is to say, we can select arbitrary $t_2$ and $t_3$, when
 \begin{equation}
     \begin{aligned}
       &  0<\alpha<\bar{\alpha}=\min \{\Theta_1,\Theta_2,\Theta_3,\Theta_4,\Theta_5,\frac{1}{L_1}  \}
     \end{aligned}
 \end{equation}
 and 
 \begin{equation}
     \begin{aligned}
         & 0<\gamma<\bar{\gamma}=\min \{\Gamma_1,\Gamma_2,\Gamma_3,\Gamma_4,\frac{1}{L_2},\frac{1}{L_3}\},
     \end{aligned}
 \end{equation}
 where $t_2>0$, $t_3>0$, $t_1>\frac{L_1t_3+L_3t_4}{\mu}$ and $t_4>\frac{2L_2t_3}{1-\rho}$, we have $\rho(R)<1$. 

Next, we use Jury criterion to seek precise range of $\alpha$ and $\beta$ to meet $\rho(Q)<1$. By computing, we can obtain the characteristic polynomial of $Q$:
 \begin{equation}
     G(\lambda)=|\lambda I-Q|=a_0+a_1\lambda+a_2\lambda^2+a_3\lambda^3+\lambda^4,
 \end{equation}
 where
 \begin{equation}
     \begin{aligned}
         a_0&=(e_1-1)[\alpha L_3(\rho\alpha L_1-e_2)+\rho^2e_3]+\rho^2\alpha\gamma e_1L_1(1+L_3),\\
         a_1&=(e_1-1)[L_3e_2(1+r)-e_4\gamma+\rho(2e_3+\rho-\alpha L_3^2)]\\
         &\quad+\gamma L_3[e_2+\alpha L_3(e_1+\rho e_1-1-2\rho)]\\
         &\quad -\rho^2e_3-\alpha L_1(1+L_3)\rho (2\gamma e_1+\rho),\\
         a_2&=(1+\gamma)L_3e_2+(e_1-1)[(1+\gamma)e_4+2\rho+e_3]\\
         &\quad+\alpha L_1[(1+L_3)(2\rho+\gamma e_1)+L_3\gamma+L_3(1+\gamma)(e_1-1-\rho)]\\
         &\quad -\gamma e_4+\rho(2e_3+\rho),\\
         a_3&=e_1-e_2-1-2\rho-(1+\gamma)e_4-\alpha L_1[L_3(2+\gamma)+1] ,
         \end{aligned}
 \end{equation}
 and $e_1=\alpha[m+L_1(1+L_3)]$, $e_2=\alpha L_2[\rho(1+L_3)-2L_3]$, $e_3=\gamma(1+\alpha L_1+\alpha L_1L_3)$ and $e_4=\alpha L_2L_3(1+L_3)$. Then we can obtain:
 \begin{equation}
     \begin{aligned}
         &b_0=a_0^2-1, b_1=a_0a_1-a_3, b_2=a_0a_2-a_2, b_3=a_0a_3-a_1\\
         &c_0=b_0^2-b_3^2,c_1=b_0b_1-b_2b_3,c_2=b_0b_2-b_1b_3.
     \end{aligned}
 \end{equation}
 Next we denote
 \begin{equation}
     \begin{aligned}
         &\mathcal{T}_1=\{(\alpha,\gamma)|a_0+a_1+a_2+a_3+1>0\},\\
        &\mathcal{T}_2=\{(\alpha,\gamma)|a_0-a_1+a_2-a_3+1>0\},\\
        &\mathcal{T}_3=\{(\alpha,\gamma)||a_0|<1\},\quad\mathcal{T}_4=\{(\alpha,\gamma)||b_0|>|b_3|\}\\
        &\mathcal{T}_5=\{(\alpha,\gamma)||c_0|>|c_2|\},\quad\mathcal{T}_6=\{(\alpha,\gamma)|\alpha>0,\gamma>0\}.
     \end{aligned}
 \end{equation}
 Thus according to Lemma $4$, when $(\alpha,\gamma) \in \bigcap\limits _{i=1}^6\mathcal{T}_i$, the spectral radius of the matrix Q is less than 1. In view of the above analysis, we know $\bigcap\limits _{i=1}^6\mathcal{T}_i$ is non-empty. Finally, denote $\rho_2=\rho(Q)$ and $0<\rho_2<1$ then we can obtain 
 \begin{equation}
     ||V_{k+1}||\leq \rho_2||V_k||.
 \end{equation}
 Furthermore, it  leads to
 \begin{equation}
     ||x_k-x^*||\leq ||V_k|| \leq C_2\rho_2^k,
 \end{equation}
 where $C_2=||V_0||$. Thus DAGT-NES can achieve the
$\mathbf{R}$-linear convergence rate. Then the proof is completed.
\end{proof}

\textbf{Remark $2$:} In Theorem $2$, we have established an $\mathbf{R}$-linear rate of DAGT-NES when the  step-size ${\alpha}$, and the largest momentum
parameter ${\gamma}$ follow  (68). Similar to DAGT-HB, the theoretical bounds of ${\alpha}$ and ${\gamma}$ in Theorem 2 are conservative. How to obtain theoretical boundaries and even optimal parameters will be considered in our future work.

\textbf{Corollary 2:} Under the same assumptions of
Theorem 2, the following equality holds:
\begin{equation}
    |f(x_k)-f(x^*)|\leq \frac{L_1}{2}C_2^2 \rho_2^{2k}.
\end{equation}
\begin{proof}Because $f(x)$ is $L_1$-smooth and $\nabla f(x^*)=0$, we can obtain
\begin{equation}
    f(x_k)-f(x^*)\leq \frac{L_1}{2}||x_k-x^*||^2.
\end{equation}
By substituting (88) into (90), we complete the proof.
\end{proof}

As a summary of this section, the DAGT-NES method is formulated as the following Algorithm \ref{algorithm2}.

\begin{algorithm}
    \caption{DAGT-NES}\label{algorithm2}
    \SetKwInOut{Input}{Input}\SetKwInOut{Output}{Output}
    \Input{initial point $x_{i,0} $ and $y_{i,0} \in \mathbb{R}^{n_i}$, $u_{i,0}=\phi_i(y_{i,0})$ and $s_{i,0}=\nabla_2 f_i(y_{i,0},u_{i,0})$ for $i=1,2,\ldots,N$, let $  (\alpha,\gamma) \in \bigcap\limits _{i=1}^6\mathcal{T}_i$, set $k=0$, select appropriate $\epsilon>0$ and $K_{max} \in \mathbb{N}^+$ .}
    \Output{optimal $x^*$ and $f(x^*)$.}
    While
    {$k<K_{max}$ and $||\nabla f(x_k)||\geq \epsilon$ do}\\
    Iterate: Update for each $i \in \{1, 2, \ldots, N\}$:
    \begin{align*}
        x_{i,k+1}\;\text{iterated by formula (58),}\\
        y_{i,k+1}\;\text{iterated by formula (59),}\\
        u_{i,k+1}\;\text{iterated by formula (60),}\\
        s_{i,k+1}\;\text{iterated by formula (61).}
    \end{align*}
    
    Update: $k=k+1$.\\
\end{algorithm}

\section{Numerical Simulation}
In order to verify the effectiveness of our proposed methods  DAGT-HB and DAGT-NES, we perform the following optimal placement as a numerical simulation. In an optimal placement problem in $\mathbf{R}^2$, suppose that there are 5 entities which are located at $r_1=(10,4)$, $r_2=(1,3)$, $r_3=(2,7)$, $r_4=(8,10)$ and $r_5=(3,9)$. And there are 5 free entities, each of which can privately be accessible to some of the fixed 5 entities. The purpose is to determine the optimal position $x_i$, $i\in \{1,2,3,4,5\}$  of the free entity so as to minimize the sum of all distances from the current position of each free entity to the corresponding fixed entity location and the distances from each entity to the weighted center of all free entities. Therefore, the cost function of  each free entity can be  modeled as follows:
\begin{equation}
    f_i(x_i,u(x))=\omega_i ||x_i-r_i||^2+ ||x_i-u(x)||^2, \quad i=1,\ldots,5,
\end{equation}
where $\omega_i$ represents the weight and is set to $20$.  We set $u(x)=\frac{\sum_{i=1}^5x_i}{5}$. So in this condition, $\phi_i$ is the identity mapping for $i=1,\ldots,5$. The communication
graph is randomly chosen to be strongly connected and doubly stochastic. 

We select the initial point $x_{1,-1}=(0,11)$, $x_{2,-1}=(9,8)$, $x_{3,-1}=(9,1)$, $x_{4,-1}=(1,4)$ and $x_{5,-1}=(3,1)$; $x_{1,0}=y_{1,0}=(2,9)$, $x_{2,0}=y_{2,0}=(8,6)$, $x_{3,0}=y_{3,0}=(7,3)$, $x_{4,0}=y_{4,0}=(4,7)$ and $x_{5,0}=y_{5,0}=(8,3)$. We initialize $u_{i,0}=\phi_i(x_{i,0})$ and $s_{i,0}=\nabla_2 f_i(x_{i,0},u_{i,0})$ for $i=1,2,\ldots,5$ and set the step size $\alpha=0.005$ and choose the momentum term $\beta=0.28$ and $\gamma=0.25$ respectively.

Then we use the DAGT-HB method to solve the optimal placement problem and the results are shown in Fig. 1 and Fig. 2. In Fig. 1,  we can see that all agents can converge to their best positions and the optimal positions are $x_1^*=(9.7524,4.1248)$, $x_2^*=(1.810,3.1714)$, $x_3^*=(2.1333,6.9810)$, $x_4^*=(7.8416,9.8381)$ and $x_5^*=(3.0857,8.8857)$ respectively.
Fig. 2  shows the  evolution of $u_{i,k}$, indicating that  the estimate $u_{i,k}$ of each free entity converges to optimal aggregative position $u(x^*)=(4.8,6.6)$ with a rapid speed.  Furthermore, Fig. 1 and Fig. 2 imply that the convergence rate is fast to support our theoretical analysis.  Meanwhile, the evolutions of $x_{i,k}$ and $u_{i,k}$ for DAGT-NES are similar to those for DAGT-HB and the results are shown in Fig. 3 and Fig. 4.  

In order to demonstrate the superiority of DAGT-HB and DAGT-NES, we compare them with DAGT under the same initial conditions. We take the first state $x_1$ as an example, the errors of  the first state are shown in Fig. 5. We find that the convergence speeds of DAGT-HB and DAGT-NES are significantly faster than that of DAGT, indicating that the introduced momentum term $x_{i,k}-x_{i,k-1}$ can enhance the convergence rate of the algorithm. Furthermore, to compare the convergence with DAGT, the losses of $(f(x_k)-f(x_*))^2$ for DAGT, DAGT-HB and DAGT-NES are shown in Fig. 6. This result further demonstrates the superiority of our algorithms. Meanwhile, we find that  DAGT-HB has a faster convergence rate but may cause more oscillations in Fig. 5. The reason of the oscillation and how to reduce it will be worth investigating in our future research.









\begin{figure}
    \centering
    \includegraphics[width=\linewidth]{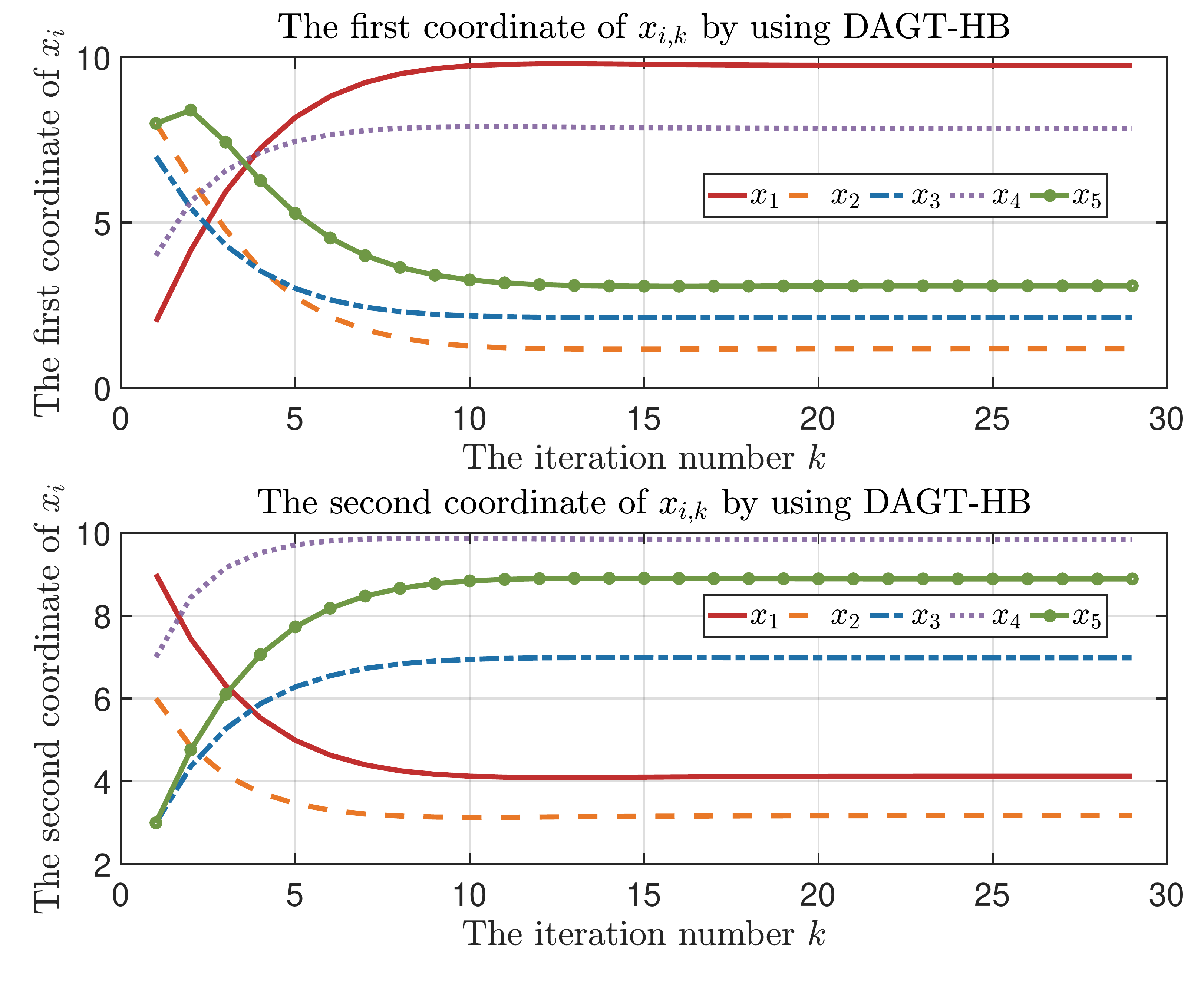}
    \caption{The evolution of $x_{i,k}$ by using DAGT-HB.}
    \label{fig:my_label}
\end{figure}

\begin{figure}
    \centering
    \includegraphics[width=\linewidth]{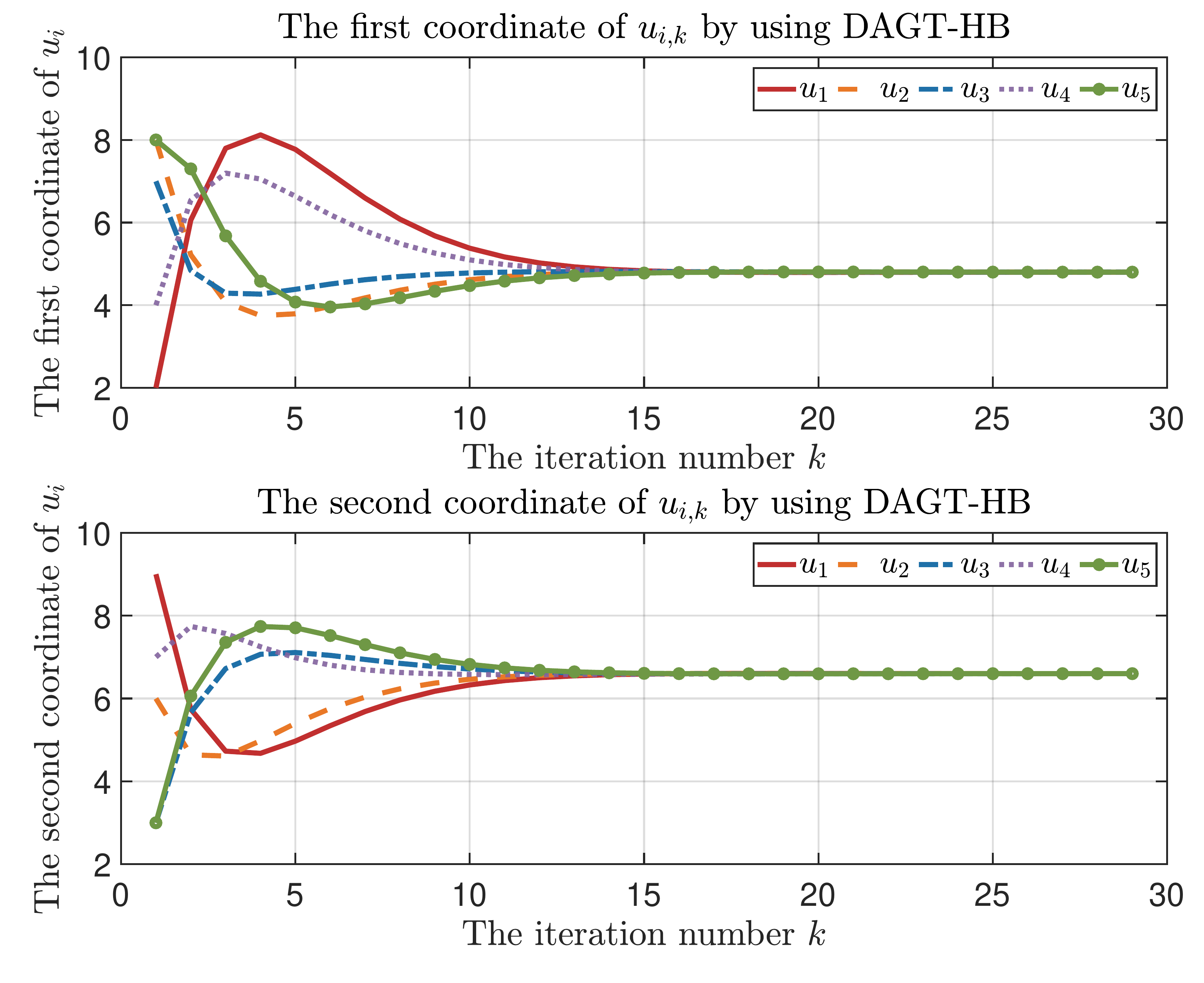}
    \caption{The evolution of $u_{i,k}$ by using DAGT-HB.}
    \label{fig:my_label}
\end{figure}

\begin{figure}
    \centering
    \includegraphics[width=\linewidth]{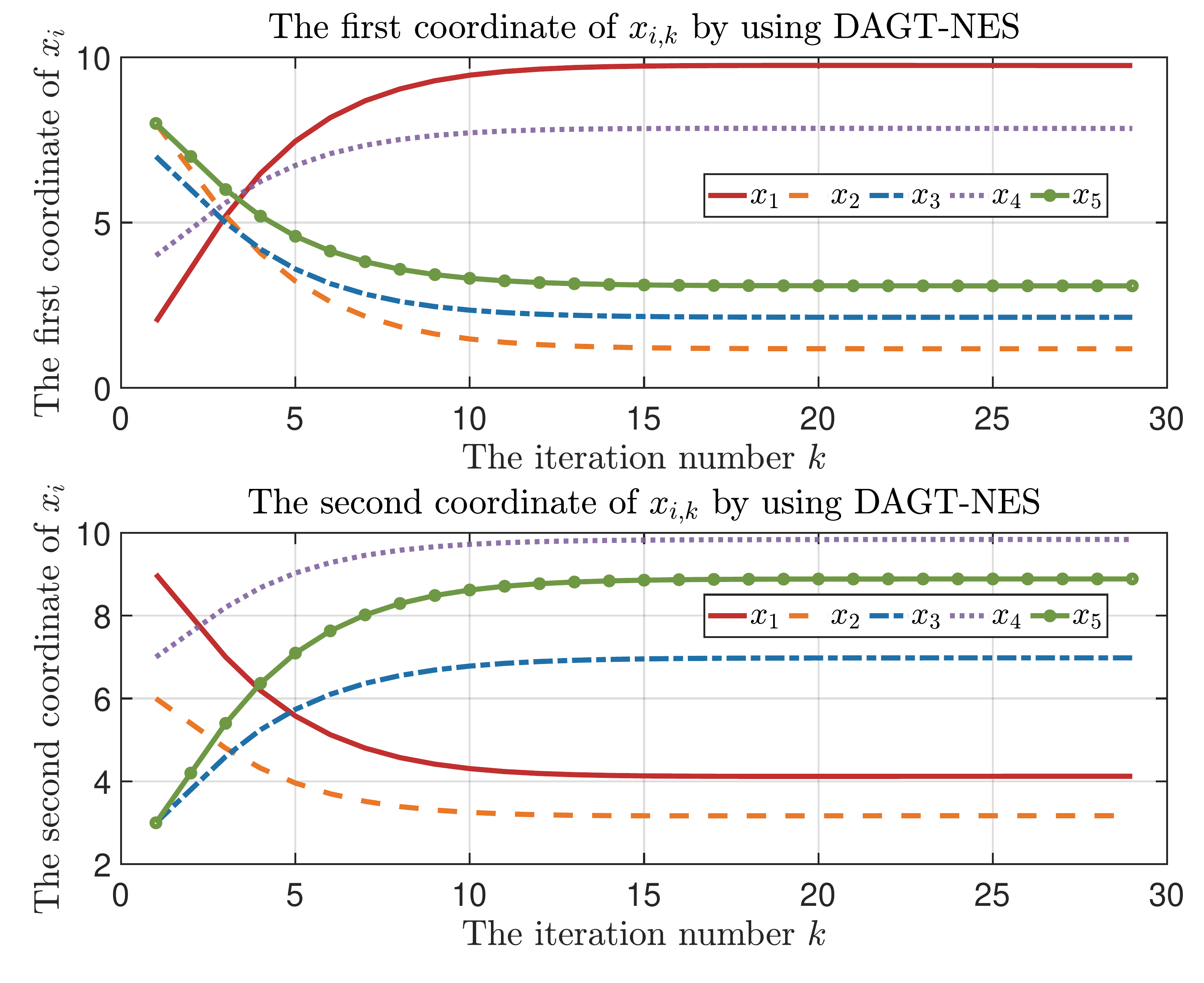}
    \caption{The evolution of $x_{i,k}$ by using DAGT-NES.}
    \label{fig:my_label}
\end{figure}

\begin{figure}
    \centering
    \includegraphics[width=\linewidth]{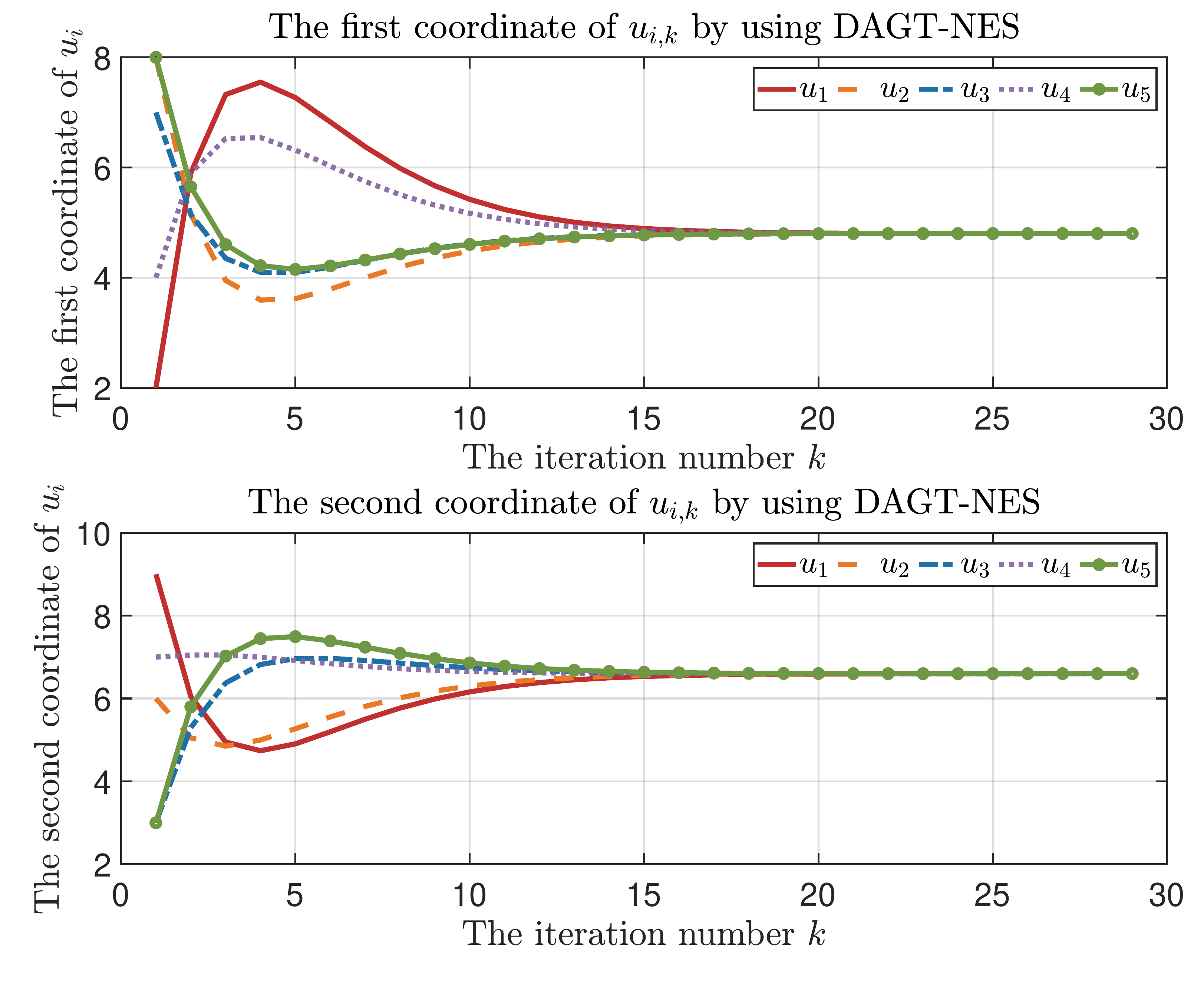}
    \caption{The evolution of $u_{i,k}$ by using DAGT-NES.}
    \label{fig:my_label}
\end{figure}

\begin{figure}
    \centering
    \includegraphics[width=\linewidth]{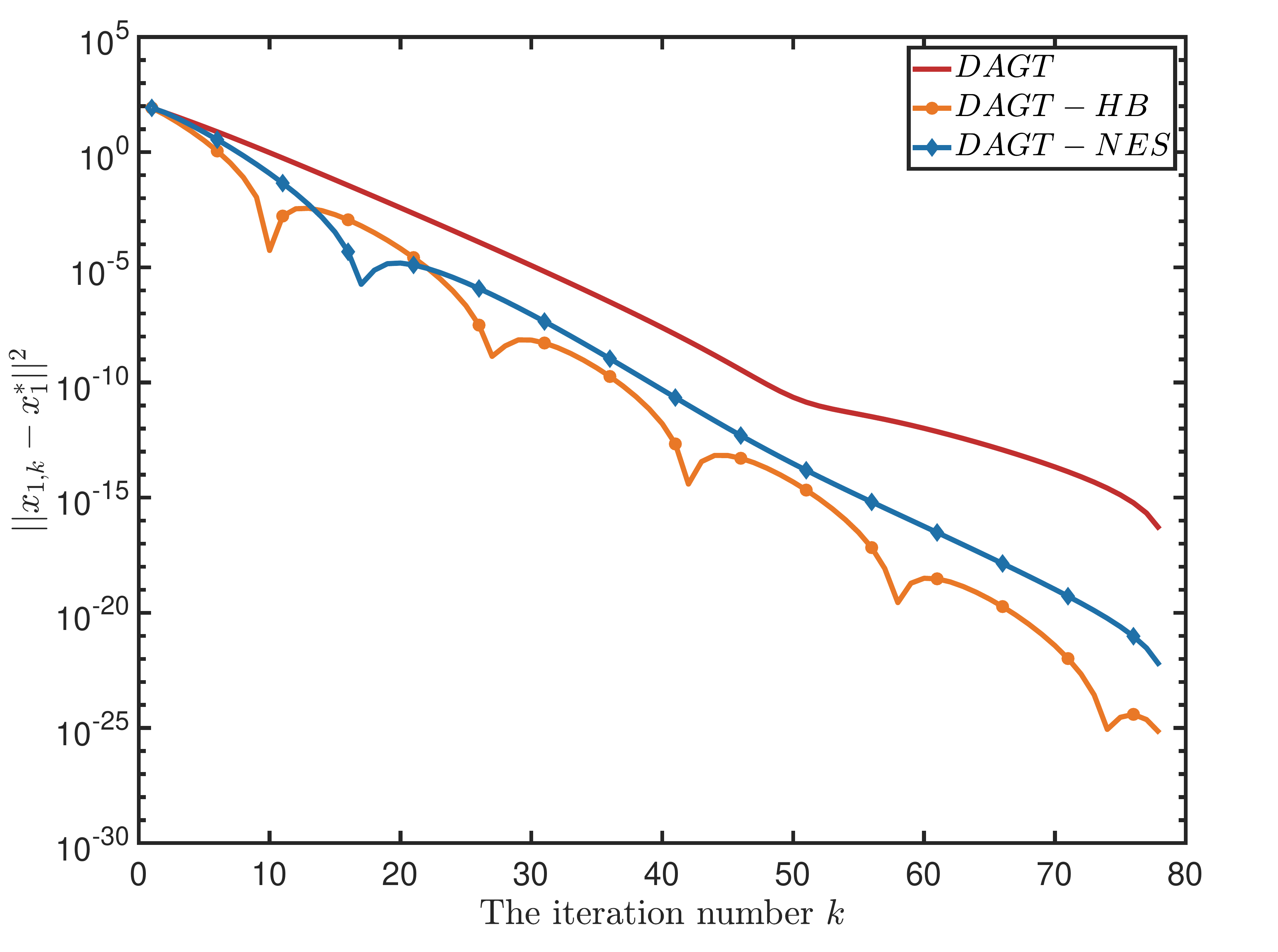}
    \caption{The state error $||x_{1,k}-x_{1}^*||^2$ comparsion among DAGT, DAGT-HB and DAGT-NES.}
    \label{fig:my_label}
\end{figure}

\begin{figure}
    \centering
    \includegraphics[width=\linewidth]{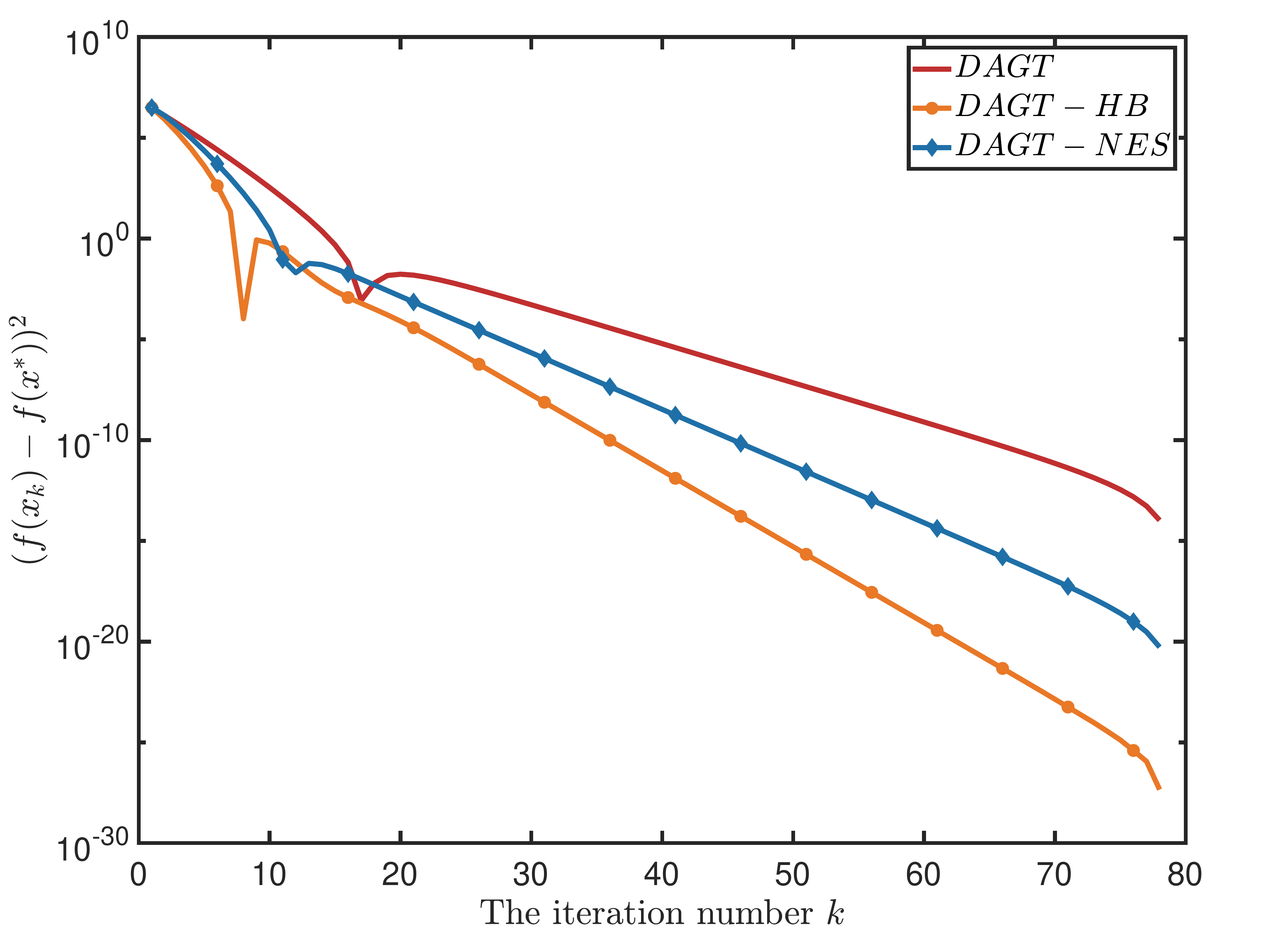}
    \caption{The cost function error $||f(x_k)-f(x^*)||^2$ comparison among DAGT, DAGT-HB and DAGT-NES.}
    \label{fig:my_label}
\end{figure}

\section{Conclusion}
This paper proposes two novel  algorithms called DAGT-HB and DAGT-NES to solve the  distributed aggregative
optimization problem in a network.   Inspired by the accelerated algorithms, we combine heavy ball and Nesterov’s accelerated method with distributed
aggregative gradient tracking method. Furthermore, we show that the algorithms DAGT-HB and DAGT-NES can converge
to an optimal solution at a global $\mathbf{R}-$linear convergence rate when
the objective function is smooth and strongly convex and when the step size and momentum term can be selected appropriately. 
Finally, we use an optimal placement problem as an example
to verify the effectiveness and
superiority of DAGT-HB and DAGT-NES. Under the same conditions, DAGT-HB and DAGT-NES can achieve much faster convergence of the state error and cost function compared to the vanilla DAGT method. 
Moreover, we find that DAGT-HB has a   faster convergence rate than DAGT-NES but may cause more oscillations. Future study may focus on the reason of oscillation, the sensitivity of parameters selection in DAGT-HB and DAGT-NES and the extension of DAGT-HB and DAGT-NES  to unbalanced graph, non-convex objective function and constrained distributed aggregation optimization problems.

\section{APPENDIX}

\subsection{Proof of Lemma 10}
Denote the equilibrium point of (58)-(61) as $x^*=$col $\left(x_1^*, \ldots, x_N^*\right)$, $y^*=$col $\left(y_1^*, \ldots, y_N^*\right)$, $u^*=$ col$\left(u_1^*, \ldots, u_N^*\right)$, and $s^*=\operatorname{col}\left(s_1^*, \ldots, s_N^*\right)$. By (59), we know $y^*=x^*$. Next, from (58), (60) and (61), we can derive
\begin{align}
    & \nabla_1 f\left(y^*, u^*\right)+\nabla \phi\left(y^*\right) s^*=\mathbf{0}_{Nd}, \\
& \mathcal{L} u^*=\mathbf{0}_{Nd},\quad \mathcal{L} s^*=\mathbf{0}_{Nd} ,
\end{align}
where $\mathcal{L}=L\otimes I_d$ and $L$ is Laplacian matrix of graph $\mathcal{G}$. Because of the property of the Laplace matrix, so $u_i^*$ must be equal to $u_j^*$ and $s_i^*$ must be equal to $s_j^*$ for all $i\neq j$. Due to equality  (62) and (63), it leads to
\begin{align}
    u_i^* & =\frac{1}{N} \sum_{i=1}^N \phi_i\left(y_{i}^*\right)=u(y^*), \\
s_i^* & =\frac{1}{N} \sum_{i=1}^N \nabla_2 f_i\left(y_{i}^*, u(y^*)\right)  .
\end{align}
By inserting (94) and (95) into (92), we can obtain $\nabla f(x^*)=\nabla f(y^*)=0$. Because $f$ is $m-$ strongly convex, $x^*$ is the unique optimal solution to problem (1).

\subsection{Proof of Lemma 11}
For $\left\|x_{k+1}-x^*\right\|$, by invoking (58), it leads to
\begin{equation}
    \begin{aligned}
        & \left\|x_{k+1}-x^*\right\| \\
& =\left\|y_k-y^*-\alpha\left[\nabla_1 f\left(y_k, u_k\right)+\nabla \phi\left(y_k\right) s_k\right]\right\| \\
& \leq \Big\Vert y_k-y^*-\alpha\left[\nabla_1 f\left(y_k, \mathbf{1}_N \otimes \bar{u}_k\right)\right. \\
&\quad \left.+\nabla \phi\left(y_k\right) \mathbf{1}_N \otimes \frac{1}{N} \sum_{i=1}^N \nabla_2 f_i\left(y_{i, k}, \mathbf{1}_N \otimes \bar{u}_k\right)\right]+\alpha \nabla f\left(y^*\right) \Big\Vert \\
&\quad +\alpha \| \nabla_1 f\left(y_k, u_k\right)+\nabla \phi\left(y_k\right) \mathbf{1}_N \otimes \bar{s}_k-\nabla_1 f\left(y_k, \mathbf{1}_N \otimes \bar{u}_k\right) \\
&\quad -\nabla \phi\left(y_k\right) \mathbf{1}_N \otimes \frac{1}{N} \sum_{i=1}^N \nabla_2 f_i\left(y_{i, k}, \mathbf{1}_N \otimes \bar{u}_k\right) \| \\
&\quad +\alpha\left\|\nabla \phi\left(y_k\right) s_k-\nabla \phi\left(y_k\right) \mathbf{1}_N \otimes \bar{s}_k\right\|. \\
    \end{aligned}
\end{equation}
From Lemma $1$, we can bound the first term of the right term of (96) as follows:
\begin{equation}
    \begin{aligned}
        &  \Big\Vert y_k-y^*-\alpha\left[\nabla_1 f\left(y_k, \mathbf{1}_N \otimes \bar{u}_k\right)\right. \\
& \quad\left.+\nabla \phi\left(y_k\right) \mathbf{1}_N \otimes \frac{1}{N} \sum_{i=1}^N \nabla_2 f_i\left(y_{i, k}, \mathbf{1}_N \otimes \bar{u}_k\right)\right]+\alpha \nabla f\left(y^*\right) \Big\Vert \\
&  \leq(1-m \alpha)\left\|y_k-y^*\right\|.
    \end{aligned}
\end{equation}
Then by inserting (59) into (97) and noting $x^*=y^*$  we can obtain
\begin{equation}
 \begin{aligned}
     \left\|y_k-y^*\right\|\leq \left\|x_k-x^*\right\|+ \gamma\left\|x_k-x_{k-1}\right\|.
 \end{aligned}   
\end{equation}
For the second term, since $f(x)$ is $L_1-$smooth and $\mathbf{1}_N \otimes \bar{u}_k=\mathcal{K} u_k$ we can get
\begin{equation}
    \begin{aligned}
        & \alpha \| \nabla_1 f\left(y_k, u_k\right)+\nabla \phi\left(y_k\right) \mathbf{1}_N \otimes \bar{s}_k-\nabla_1 f\left(y_k, \mathbf{1}_N \otimes \bar{u}_k\right) \\
&\quad -\nabla \phi\left(y_k\right) \mathbf{1}_N \otimes \frac{1}{N} \sum_{i=1}^N \nabla_2 f_i\left(y_{i, k}, \mathbf{1}_N \otimes \bar{u}_k\right) \| \\
& \leq \alpha L_1\left\|u_k-\mathcal{K} u_k\right\|.\\
    \end{aligned}
\end{equation}
For the last term, by using Assumption $4$ and $\mathbf{1}_N \otimes \bar{s}_k=\mathcal{K} s_k$ we can  the following inequality:
\begin{equation}
    \alpha\left\|\nabla \phi\left(y_k\right) s_k-\nabla \phi\left(y_k\right) \mathbf{1}_N \otimes \bar{s}_k\right\| \leq \alpha L_3\left\|s_k-\mathcal{K} s_k\right\|.
\end{equation}
Then by bonding (97)-(100) with (96), we complete the proof.

\subsection{Proof of Lemma 12}
For $\left\|x_{k+1}-x_k\right\|$, we invoke (58) and note  $\nabla f(x^*)=0$, then we have
\begin{equation}
    \begin{aligned}
        & \left\|x_{k+1}-x_k\right\| \\
& =\left\|\gamma(x_k-x_{k-1})-\alpha(\nabla_1 f\left(y_k, u_k\right)+\nabla \phi\left(y_k\right) s_k)\right\| \\
& \leq \alpha \Big\Vert \nabla_1 f\left(y_k, u_k\right)+\nabla \phi\left(y_k\right) \mathcal{K} s_k-\nabla_1 f\left(y^*, \mathbf{1}_N \otimes u^*\right) \\
& \quad\quad-\nabla \phi\left(x^*\right)\left[\mathbf{1}_N \otimes \frac{1}{N} \sum_{i=1}^N \nabla_2 f_i\left(y^*, \mathbf{1}_N \otimes u^*\right)\right] \Big\Vert \\
& \quad+\alpha\left\|\nabla \phi\left(y_k\right)\left(s_k-\mathcal{K} s_k\right)\right\|+\gamma\left\|x_k-x_{k-1}\right\|. \\
    \end{aligned}
\end{equation}
By utilizing Assumption $2$ and triangle inequality of norm, we can obtain the following formula:
\begin{equation}
    \begin{aligned}
        &   \Big\Vert \nabla_1 f\left(y_k, u_k\right)+\nabla \phi\left(y_k\right) \mathcal{K} s_k-\nabla_1 f\left(y^*, \mathbf{1}_N \otimes u^*\right) \\
& \quad-\nabla \phi\left(y^*\right)\left[\mathbf{1}_N \otimes \frac{1}{N} \sum_{i=1}^N \nabla_2 f_i\left(y^*, \mathbf{1}_N \otimes u^*\right)\right] \Big\Vert \\
& \leq  L_1\left(\left\|y_k-x^*\right\|+\left\|u_k-\mathbf{1}_N \otimes u^*\right\|\right)\\
& \leq  L_1\left(\left\|x_k-x^*\right\|+\left\|u_k-\mathcal{K} u_k\right\|\right) \\
& \quad+\gamma L_1\left\|x_k-x_{k-1}\right\|+L_1\left\|\mathcal{K} u_k-\mathbf{1}_N \otimes u^*\right\|.
    \end{aligned}
\end{equation}
By (30), we know
\begin{equation}
    \begin{aligned}
        \left\|\mathcal{K} u_k-1_N \otimes u^*\right\| & \leq L_3\left\|y_k-y^*\right\|\\
        & \leq  L_3\left\|x_k-x^*\right\|+ L_3\gamma\left\|x_k-x_{k-1}\right\|.
    \end{aligned}
\end{equation}
Then by using Assumption $4$ we can obtain
\begin{equation}
    \left\|\nabla \phi\left(y_k\right)\left(s_k-\mathcal{K} s_k\right)\right\| \leq L_3 \left\|s_k-\mathcal{K} s_k\right\|.
\end{equation}
Finally, by inserting (102)-(104)  into (101) then we  can obtain the Lemma $12$.

\subsection{Proof of Lemma 13}
For $\left\|u_{k+1}-\mathcal{K} u_{k+1}\right\|$, by invoking (60), it leads to
\begin{equation}
\begin{aligned}
& \left\|u_{k+1}-\mathcal{K} u_{k+1}\right\| \\
& =\left\|\mathcal{A} u_k+\phi\left(y_{k+1}\right)-\phi\left(y_k\right)-\mathcal{K} \mathcal{A} u_k-\mathcal{K}\left[\phi\left(y_{k+1}\right)-\phi\left(y_k\right)\right]\right\| \\
& \leq \rho\left\|u_k-\mathcal{K} u_k\right\|+\|I-\mathcal{K}\|\left\|\phi\left(y_{k+1}\right)-\phi\left(y_k\right)\right\| \\
& \leq \rho\left\|u_k-\mathcal{K} u_k\right\|+L_3\left\|y_{k+1}-y_k\right\|\\
& = \rho\left\|u_k-\mathcal{K} u_k\right\|+L_3\left\|(1+\gamma)(x_{k+1}-x_k)-\gamma (x_k-x_{k-1})\right\|\\
& \leq\rho\left\|u_k-\mathcal{K} u_k\right\|+L_3(1+\gamma)\left\|x_{k+1}-x_k\right\|+\gamma L_3\left\|x_k-x_{k-1}\right\|,\\
\end{aligned}
\end{equation}
where Lemma $2$ has been utilized to obtain the first inequality, and by using  Assumption $4$ we can obtain the third inequality. Then by substituting (65) into (105) we can  obtain the Lemma $13$.

\subsection{Proof of Lemma 14}
 For $\left\|s_{k+1}-\mathcal{K} s_{k+1}\right\|$, by invoking (61)  we can obtain
\begin{equation}
\begin{aligned}
& \left\|s_{k+1}-\mathcal{K} s_{k+1}\right\| \\
& =||\mathcal{A} s_k+\nabla_2 f\left(y_{k+1}, u_{k+1}\right)-\nabla_2 f\left(y_k, u_k\right)\\
&\quad-\mathcal{K}\mathcal{A} s_k-\mathcal{K}[\nabla_2 f\left(y_{k+1}, u_{k+1}\right)-\nabla_2 f\left(y_k, u_k\right)]|| \\
& \leq \left\|\mathcal{A}s_k-\mathcal{K} s_k\right\|\\
& \quad+\|I-\mathcal{K}\|\left\|\nabla_2 f\left(y_{k+1}, u_{k+1}\right)-\nabla_2 f\left(y_k, u_k\right)\right\|\\
& \leq \rho\left\|s_k-\mathcal{K} s_k\right\|+\left\|\nabla_2 f\left(y_{k+1}, u_{k+1}\right)-\nabla_2 f\left(y_k, u_k\right)\right\| \\
& \leq \rho\left\|s_k-\mathcal{K} s_k\right\|+  L_2\left(\left\|y_{k+1}-y_k\right\|+\left\|u_{k+1}-u_k\right\|\right),\\
\end{aligned}
\end{equation}
where Assumption $2$ has been leveraged in the last inequality. Then by (36), we know
\begin{equation}
\begin{aligned}
    & \left\|u_{k+1}-u_k\right\|
 \leq 2 \left\|u_k-\mathcal{K} u_k\right\|+ L_3\left\|y_{k+1}-y_k\right\|.
    \end{aligned}
\end{equation}
By substituting  (107) into (106), it can lead to
\begin{equation}
    \begin{aligned}
        & \left\|s_{k+1}-\mathcal{K} s_{k+1}\right\| \\
        & < \rho\left\|s_k-\mathcal{K} s_k\right\|+  L_2(L_3+1)\left\|y_{k+1}-y_k\right\|+2L_2\left\|u_k-\mathcal{K} u_k\right\| \\
        & \leq \rho\left\|s_k-\mathcal{K} s_k\right\|+ L_2(L_3+1)(\gamma+1)\left\|x_{k+1}-x_k\right\|\\
        & \quad+ L_2(L_3+1)\gamma\left\|x_{k}-x_{k-1}\right\|+2L_2\left\|u_k-\mathcal{K} u_k\right\|.
    \end{aligned}
\end{equation}
Then by substituting (65) into (108), we  can be obtain Lemma $14$.

\bibliographystyle{unsrt}
\bibliography{reference}  

\end{document}